\documentclass[12pt,a4paper]{article}
\usepackage{amsfonts,amsmath,amssymb,amscd}
\usepackage{fancybox}
\usepackage[latin1]{inputenc}
\usepackage[usenames]{color}
\usepackage{graphicx}
\usepackage{float}
\usepackage{fullpage}
\usepackage{epsfig,amssymb}
\usepackage{enumerate}
\usepackage{subcaption}

\usepackage[normalem]{ulem}
\definecolor{avocado}{rgb}{0.34,0.51,0.01}

\newtheorem{thm}{Theorem}[section]

\newtheorem{lem}[thm]{Lemma}
\newtheorem{prop}[thm]{Proposition}

\newtheorem{rem}[thm]{Remark}

\begin{document}
	
\title{\textbf{Instability of Closed $p$-Elastic Curves in $\mathbb{S}^2$}}

\author{A. Gruber, A. P\'ampano and M. Toda}
\date{\today}

\maketitle 

\begin{abstract}
\noindent For $p\in\mathbb{R}$, we show that non-circular closed $p$-elastic curves in $\mathbb{S}^2$ exist only when $p=2$, in which case they are classical elastic curves, or when $p\in(0,1)$. In the latter case, we prove that for every pair of relatively prime natural numbers $n$ and $m$ satisfying $m<2n<\sqrt{2}\,m$, there exists a closed spherical $p$-elastic curve with non-constant curvature which winds around a pole $n$ times and closes up in $m$ periods of its curvature. Further, we show that all closed spherical $p$-elastic curves for $p\in(0,1)$ are unstable as critical points of the $p$-elastic energy.
\\

\noindent{\emph{Keywords:} closed curves, $p$-elastic curves, $p$-elastic energy, stability.}
\end{abstract}

\section{Introduction}

Functionals depending on curvature invariants occupy a central position in the study of variational problems for curves, and have been used as effective scientific models throughout their long and storied history.  

Originating from the elastic theory pioneered by the Bernoulli family and L. Euler, such models were essential in the early developments of the Calculus of Variations. The problem of quantifying the bending deformation of rods was first formulated by J. Bernoulli in 1691 \cite{L} and led to the complete description of the classical elastic curves by L. Euler.  Printed in an Appendix to his famous 1744 monograph \cite{E}, these curves represent the possible qualitative types for rods in untwisted planar configurations, unifying and extending the special cases known to J. Bernoulli around 1694 \cite{Be}.

To accomplish this goal, L. Euler followed a general principle formulated originally by D. Bernoulli: an elastic rod should bend along the curve which minimizes the potential energy of the strain under suitable constraints. Consequently, he suggested finding extrema of the total squared curvature functional. However, this idea can easily be posed more generally: in a letter to L. Euler of 1738 \cite{T}, D. Bernoulli proposed to investigate extrema of the functionals (in modern notation)
$$\mathbf{\Theta}_p(\gamma):=\int_\gamma \kappa^p\,ds\,,$$
where $\kappa$ is the curvature of the curve $\gamma$ and $s$ is its arc length parameter.  This gives a rich family of curves which has yet to be fully understood.

Clearly, $p=2$ recovers the classical bending energy whose extrema are (free) elastic curves. Since their discovery, these curves have received a considerable amount of attention because they arise in a vast number of applications. They are related to the classical \emph{lintearia} which represents the shape of a long cloth sheet full of water \cite{L}, they are the generating curves of the first non-conformally minimal Willmore tori \cite{Pi}, they serve as generating curves for biomembranes in the Canham-Helfrich-Evans models for lipid bilayers \cite{C,Ev,H}, and they are used for smoothing and pathfinding in computer vision \cite{M}, to mention just a few. Moreover, elastic curves also arise in various mathematical contexts: they serve in integrable flows as the traveling wave solutions of both the mKdV equation \cite{GP1,GP2,LPe} and (after the Hasimoto transformation) the non-linear Schr\"{o}dinger equation \cite{Ha1,Ha2,K}, and find similar uses in elastic flows \cite{DLLPS,DKS}.

At this point, elastic curves have been studied through a variety of different techniques. Using a Lagrangian approach and the theory of Killing vector fields, the existence of closed elastic curves was studied in \cite{LS}. From a different perspective, elastic curves were analyzed following Griffith's approach to calculus of variations in \cite{BG}. We highlight here that these curves may be viewed as solutions of an optimal control problem, and so they are strongly connected to the Maximum Principle and its associated Hamiltonian (for details, see \cite{J}).

Chronologically, the second particular case to receive widespread attention is $p=0$, for which $\mathbf{\Theta}_p$ is simply the arc length functional. This variational problem was stated in 1697 as a public challenge from Johan Bernoulli to Jacob Bernoulli \cite{L}. Nowadays, it is well known that the critical curves for $\mathbf{\Theta}_0$ are geodesics, that is, those curves whose (geodesic) curvature is identically zero. This represents a far-reaching generalization of locally shortest path which reduces to straight lines in Euclidean spaces.

Moving to $p=1$ yields the total curvature type energy. Since the total curvature of closed curves in $\mathbb{R}^2$ is a topological invariant, this variational problem is trivial when the functional is acting on closed planar curves. Nevertheless, non-trivial extensions of this energy have been widely studied for applications to other fields. For instance, in visual curve completion the primary visual cortex is usually modeled as the unit tangent bundle of the plane with a suitable sub-Riemannian geometry, and it is believed that geodesics in this space are optimal for image reconstruction \cite{BYBS,Pe}. These geodesics project down to minimizers of a total curvature type energy in the plane \cite{AGP0,BYBS}.

The case $p=1/2$ was considered by W. Blaschke in 1921 \cite{B}, who explicitly obtained the curvature of critical curves in terms of their arc length parameter. He showed that planar critical curves are catenaries. In \cite{AGP}, critical curves for an extension of this functional were characterized as the profile curves of invariant constant mean curvature (CMC) surfaces in Riemannian and Lorentzian $3$-space forms. In particular, the functional $\mathbf{\Theta}_{1/2}$ arises precisely in the minimal case.

In 1923, W. Blaschke also considered minimizers of the functional $\mathbf{\Theta}_p$ for $p=1/3$ \cite{B}. This functional represents the equi-affine length for convex curves whose critical curves are parabolas. Recently, the equi-affine geometry of convex curves has been consistently used in studies on human curvilinear two-dimensional drawing movements \cite{FH} as well as recognition for non-rigid planar shapes \cite{RK}.

For integer values of $p>2$, this variational problem was considered a couple of decades ago (see e.g. \cite{AGM}) and has been useful in constructing  Willmore-Chen submanifolds in certain Riemannian and pseudo-Riemannian spaces with warped product metrics \cite{ABG,BFLM}, as well as for analyzing conformal tension in string theories \cite{BFL}.

For suitable rational values of $p\in[1/4,1)$, critical curves of $\mathbf{\Theta}_p$ can be characterized as the generating curves of non-CMC rotational bi-conservative hypersurfaces in space forms of arbitrary dimension (see \cite{MOP,MP}). Bi-conservative hypersurfaces are defined as those hypersurfaces for which the stress-energy tensor associated to the bi-energy is conservative, which can be seen as a geometric generalization of the CMC condition.

In general, for $p\neq 0,1$, the critical curves for these variational problems completely determine the rotational surfaces in $\mathbb{R}^3$ satisfying a linear relation between their principal curvatures \cite{LP2}. In fact, it was shown in \cite{PhD} that the same result holds in any Riemannian $3$-space form. Surfaces satisfying a linear relation between their principal curvatures are referred as linear Weingarten surfaces in the literature. In particular, some of these rotational linear Weingarten surfaces are generalized bi-conservative surfaces \cite{LYZ}.

The result of \cite{LP2} generalizes the well known relation between the (free) elastic curve ($p=2$) in $\mathbb{R}^2$ and the Mylar balloon, which is the physical object representing the closed rotational surface in $\mathbb{R}^3$ that maximizes volume while preserving the length of the profile curve \cite{MO}. Moreover, the same curves have just been characterized as the generating curves of cylinders in $\mathbb{R}^3$, which are critical for the area functional and a suitable vertical potential energy \cite{LP}. This last problem includes the classical \emph{lintearia} for the case where the vertical potential density is proportional to the height, i.e., under a constant gravitational field.

In the papers \cite{LP2} and \cite{LP}, the authors described the curves in $\mathbb{R}^2$ geometrically and depicted their shape. Interestingly, one can easily prove that these planar critical curves are never closed. Except for the trivial variational problem corresponding to $p=1$, critical curves in $\mathbb{R}^2$ for $\mathbf{\Theta}_p$ are generically not closed (for more details on the closure condition, one can see \cite{PhD}). Consequently, one is tempted to look at the variational problem preserving the length of these curves. There is a well known non-trivial closed critical curve when $p=2$, namely the elastic figure-eight \cite{L}, and some more general closed critical curves for $\mathbf{\Theta}_p$ with length constraint were obtained in \cite{LP} (see also \cite{MuP} for the case $p=1/2$ in the sphere).

However, in order to find closed critical curves for $\mathbf{\Theta}_p$ without length constraint, it is natural to look at the spherical case. In the round $2$-sphere $\mathbb{S}^2$ there exist closed (free) elastic curves ($p=2$). In fact, there exists a discrete bi-parametric family of closed non-circular elastic curves, including simple ones \cite{LS}. In the case $p=1/2$, the existence of countably many non-trivial closed critical curves in $\mathbb{S}^2$ was first proved in \cite{AGM}. Later on, in \cite{AGP1} using an argument involving elliptic integrals, it was shown that this family of curves is in one-to-one correspondence with pairs of relatively prime natural numbers $n$ and $m$ satisfying $m<2n<\sqrt{2}\,m$. As a consequence, none of these curves is simple. An existence result for infinitely many rational numbers $p\in[1/4,1)$ was also obtained in \cite{MOP} and \cite{MP}. Quite surprisingly, the restriction on the parameters $n$ and $m$ turns out to be the same as in the case $p=1/2$ and independent from the value of $p$. However, the existence of closed critical curves for $\mathbf{\Theta}_p$ in $\mathbb{S}^2$ cannot be taken for granted, since it is not always guaranteed. For instance, when $p>2$ is a natural number the only closed critical curves are geodesics \cite{AGM}. Note that regularity assumptions are essential in this assertion: in \cite{AGM} curves were assumed to be of class at least $\mathcal{C}^4$ (which will also be the case of the present paper). The existence of non-geodesic critical curves with less regularity have been shown recently in \cite{SW}.

Motivated by the above mentioned results, this work investigates the question of whether or not (non-trivial) closed spherical critical curves for $\mathbf{\Theta}_p$ exist for a given parameter $p\in\mathbb{R}$.  Moreover, in the event that such curves exist, we aim to classify whether or not these curves are stable. 

With these goals in mind, Section 2 recalls the Euler-Lagrange equation associated to $\mathbf{\Theta}_p$ and proves Theorem \ref{restriction}, which asserts that critical curves with non-constant \emph{periodic} curvature (a necessary, but not sufficient, condition to get closed curves) only arise when $p=2$ or $p\in(0,1)$. Aside from the well known case $p=2$ \cite{LS}, we also show in Proposition \ref{periodic} that all critical curves for $\mathbf{\Theta}_p$ with $p\in(0,1)$ have periodic curvature (assuming they are defined on their maximal domains). Consequently, these curves are the ideal candidates for which one should check the closure condition. Section 3 is devoted to showing the existence of closed curves for every $p\in(0,1)$ by analysis of the corresponding closure condition. Indeed, we prove that for every pair of relatively prime natural numbers $n$ and $m$ satisfying $m<2n<\sqrt{2}\,m$ (a condition surprisingly independent of the value of $p$) there exists a (non-trivial) closed spherical curve critical for $\mathbf{\Theta}_p$ (Theorem \ref{existence}). In fact, there is  convincing numerical evidence that, for fixed $n$ and $m$, this curve is unique up to isometries of the sphere. Moreover, the numerical experiments strongly support the \emph{ansatz} that all (non-trivial) closed critical curves for $\mathbf{\Theta}_p$ with $p\in(0,1)$ are one of the above described curves. In particular, if $p\in(0,1)$, we may conclude from these experiments that (non-trivial) simple and closed critical curves for $\mathbf{\Theta}_p$ will not exist.

In the case $p=2$, it turns out that the only stable and closed critical curve is a geodesic \cite{LS}. It is then reasonable to check if a similar result holds for $p\in(0,1)$. In Section 4 we analyze the second variation formula of $\mathbf{\Theta}_p$ to prove that, surprisingly, all closed critical curves are unstable (Theorem \ref{instability}).

\section{Critical Curves}

Let $(x,y,z)$ be the standard coordinates of the Euclidean space $\mathbb{R}^3$. The round $2$-sphere of constant sectional curvature $\rho>0$ can be understood as the hyperquadric
$$\mathbb{S}^2(\rho)=\{(x,y,z)\in\mathbb{R}^3\,\lvert\, x^2+y^2+z^2=1/\rho\}\,,$$
endowed with the induced metric from  $\mathbb{R}^3$. Without any significant loss of generality, from now on we will assume $\rho=1$ and use the notation $\mathbb{S}^2\equiv\mathbb{S}^2(1)$.

Let $\gamma:I\subseteq\mathbb{R}\longrightarrow\mathbb{S}^2$ be a smooth immersed (spherical) curve parameterized by the arc length parameter $s\in I$. Denote by $T(s):=\gamma'(s)$ the unit tangent vector field along the curve $\gamma(s)$, where $\left(\,\right)'$ denotes the derivative with respect to the arc length parameter $s$, and define the unit normal vector field $N(s)$ along $\gamma(s)$ to be the counter-clockwise rotation of $T(s)$ through an angle $\pi/2$ in the tangent bundle of $\mathbb{S}^2$. In this setting, the (signed) geodesic \emph{curvature} $\kappa(s)$ of $\gamma(s)$ is defined by the Frenet-Serret equation
$$\nabla_T T(s)=\kappa(s)\,N(s)\,,$$
where $\nabla$ denotes the Levi-Civita connection on the tangent bundle to $\mathbb{S}^2$. It follows from the Fundamental Theorem of Curves that the curvature $\kappa(s)$ completely determines the (spherical) curve $\gamma(s)$ up to isometries of $\mathbb{S}^2$.

The main object of study is the \emph{$p$-elastic functional}, $p\in\mathbb{R}$,
\begin{equation}\label{energy}
\mathbf{\Theta}_p(\gamma):=\int_\gamma \kappa^p\,ds\,,
\end{equation}
which acts on the space of smooth immersed curves $\gamma:I\subseteq\mathbb{R}\longrightarrow\mathbb{S}^2$ parameterized by the arc length $s\in I$ and denoted by $\mathcal{C}^\infty(I,\mathbb{S}^2)$.  Note that although we request that the curves under consideration be smooth, it is enough to consider only $\mathcal{C}^4$ regularity. If $p\in\mathbb{R}\setminus\mathbb{N}$ (we consider $0\in\mathbb{N}$), $\mathbf{\Theta}_p$ is understood to act on the subspace of convex curves $\mathcal{C}_*^\infty(I,\mathbb{S}^2)$, i.e., it is assumed that $\kappa>0$ holds for all curves in the variation.

From standard arguments involving integration by parts, one can compute the first variation formula associated to $\mathbf{\Theta}_p$ and conclude that a critical curve must satisfy the following Euler-Lagrange equation on its interior,
\begin{equation}\label{EL}
p\,\frac{d^2}{ds^2}\left(\kappa^{p-1}\right)+\left(p-1\right)\kappa^{p+1}+p\,\kappa^{p-1}=0\,.
\end{equation}
Those curves whose curvature function is a solution of \eqref{EL} are called \emph{$p$-elastic curves}.

If $p=0$, the Euler-Lagrange equation \eqref{EL} reduces to $\kappa=0$ and so  the only $0$-elastic curves are geodesics, which are congruent to the equator. 
As the functional $\mathbf{\Theta}_0$ is just the length functional, this conclusion is expected. On the other hand, if $p=1$ it is clear from \eqref{EL} that critical curves do not exist. Observe that this is a huge difference from the Euclidean case, where the associated Euler-Lagrange equation of $\mathbf{\Theta}_1$ is an identity.

The present study begins with the existence of $p$-elastic circles, i.e., constant solutions of \eqref{EL}.

\begin{prop} Assume that $\gamma$ is a circle critical for $\mathbf{\Theta}_p$. Then, either $p\in\mathbb{N}$ $(p\neq 1)$ and $\gamma$ is a geodesic, or $p\in[0,1)$ and the curvature of $\gamma$ is
\begin{equation}\label{circles}
\kappa=\sqrt{\frac{p}{1-p}}\,.
\end{equation}
\end{prop}
\textit{Proof.} Let $\gamma$ be a critical circle for $\mathbf{\Theta}_p$. From criticality it follows that the Euler-Lagrange equation \eqref{EL} must be satisfied. Moreover, since $\gamma$ is a circle, its curvature $\kappa$ must be constant and \eqref{EL} simplifies to
$$\kappa^{p-1}\left(\left[p-1\right]\kappa^2+p\right)=0\,.$$
From this, either $\kappa=0$ which is only possible when $p\neq 1$ is a natural number, or the curvature satisfies
$$\kappa^2=\frac{p}{1-p}\,.$$
The latter condition necessarily implies $p\in [0,1)$ in order for the right hand side to be non-negative. This gives rise to the second option in the statement. \hfill$\square$
\\

When $p\in[0,1)$, it is convenient to understand $p$-elastic circles in terms of their corresponding Euclidean radius $r$. The Gauss equation implies that the curvature $\bar{\kappa}$ of $\gamma$ as a curve in $\mathbb{R}^3$ satisfies the following relation with the (geodesic) curvature $\kappa$ in $\mathbb{S}^2$,
$$\bar{\kappa}^2=\kappa^2+1\,.$$
Consequently, from \eqref{circles} it follows that
$$\bar{\kappa}^2=\frac{p}{1-p}+1=\frac{1}{1-p}\,.$$
Using the relation $r=1/\bar{\kappa}$ for the radius of curvature of a planar curve, this implies that the Euclidean radii of $p$-elastic circles are given by
\begin{equation}\label{radius}
r=\sqrt{1-p}\,.
\end{equation}
Note that for the case $p=0$ the above radius is one, and hence coincides with the radius of the equator of $\mathbb{S}^2$, as expected. Moreover, when $p\to 1$, the radius \eqref{radius} tends to zero meaning that $p$-elastic circles converge to a point.

Of course, it is more interesting to consider the case where the curvature of $p$-elastic curves is not constant. This yields a useful conservation law for the Euler-Lagrange equation \eqref{EL}.

\begin{prop}\label{ODEprop} Let $\gamma(s)$ be a $p$-elastic curve with nonconstant curvature $\kappa(s)$. Then, $\kappa(s)$ is a solution of the first order ordinary differential equation
\begin{equation}\label{ODE}
p^2(p-1)^2\kappa^{2(p-2)}\left(\kappa'\right)^2+(p-1)^2\kappa^{2p}+p^2\kappa^{2(p-1)}=a\,,
\end{equation}
where $a>0$ is a constant of integration.
\end{prop}
\textit{Proof.} Since the curvature of $\gamma(s)$ is nonconstant, its derivative with respect to $s$ is not identically zero. Multiplying the Euler-Lagrange equation \eqref{EL} by $\kappa'$ and suitably manipulating the result produces an exact differential equation which integrates to \eqref{ODE}.  \hfill$\square$

\begin{rem} An alternative proof of Proposition \ref{ODEprop} can be obtained by adapting the computations of \cite{LS} involving Killing vector fields along curves (see also \cite{AGM,PhD} for more details). Following the Hamiltonian approach, the constant of integration $a>0$ has a clear physical meaning: it represents the square of the length of the momentum of the critical curve. The momentum is constant (conserved) along $p$-elastic curves due to Proposition \ref{ODEprop}, and can be identified with the sum of the (adapted) Killing vector fields along $\gamma$.
\end{rem}

Since the first goal of the paper is to prove the existence of closed $p$-elastic curves, the main emphasis will be placed on periodic solutions to \eqref{ODE}. However, the existence of periodic solutions to the ordinary differential equation \eqref{ODE} is highly restricted. In fact, we have already seen that for $p=0,1$, there are no critical curves with non-constant curvature for $\mathbf{\Theta}_p$, and so there are no solutions (periodic or not) to \eqref{ODE}. Moreover, in \cite{AGM} it was shown that if $p>2$ is a natural number, then the only $p$-elastic curves with periodic curvature are geodesics. The following proposition extends these results to arbitrary values of $p\in\mathbb{R}$.

\begin{thm}\label{restriction} Let $\gamma$ be a $p$-elastic curve with non-constant periodic curvature. Then, either $p=2$ or $p\in(0,1)$.
\end{thm}
\textit{Proof.} Consider a $p$-elastic curve $\gamma$ with non-constant periodic curvature. The cases $p=0,1$ have already been excluded from the previous discussion. Similarly, the non-existence of $p$-elastic curves for the case when $p>2$ is a natural number was shown in Proposition 8 of \cite{AGM} using arguments about the positive roots of a polynomial, $\widetilde{Q}_{p,a}(\kappa)$ (see below). The case $p=2$ is well known and existence was shown in \cite{LS}. For the other cases, i.e., $p\in\mathbb{R}\setminus\mathbb{N}$, recall that these curves are restricted to being convex. In other words, the curvatures of these $p$-elastic curves must satisfy $\kappa>0$. 

Assume first that $p<0$ holds, so that \eqref{ODE} can be rewritten as
$$\left(\kappa'\right)^2=\frac{\kappa^2}{p^2(1-p)^2}\left(a\,\kappa^{2(1-p)}-[1-p]^2\kappa^2-p^2\right)=\frac{\kappa^2}{p^2(1-p)^2}Q_{p,a}(\kappa)\,.$$
In order for periodic solutions to exist, the equation $Q_{p,a}(\kappa)=0$ must have two positive solutions. However a simple differentiation shows that $Q_{p,a}(\kappa)$ has just one critical point for $\kappa>0$, which is a local minimum. In combination with $Q_{p,a}(0)=-p^2<0$, this proves that $Q_{p,a}(\kappa)=0$ cannot have two positive solutions.

Assume now that $p\in(1,2)$. As before, \eqref{ODE} rewrites as
$$\left(\kappa'\right)^2=\frac{\kappa^{2(2-p)}}{p^2(p-1)^2}\left(a-\left[p-1\right]^2\kappa^{2p}-p^2\kappa^{2(p-1)}\right)=\frac{\kappa^{2(2-p)}}{p^2(p-1)^2}\widetilde{Q}_{p,a}(\kappa)\,,$$
and again $\widetilde{Q}_{p,a}(\kappa)=0$ cannot have two positive solutions. In this case, this assertion follows directly from the monotonicity of $\widetilde{Q}_{p,a}(\kappa)$ for $\kappa>0$.

Finally, consider the case $p\in\mathbb{R}\setminus\mathbb{N}$ and $p>2$. Then, \eqref{ODE} can be written as
$$p^2(p-1)^2\kappa^{2(p-2)}\left(\kappa'\right)^2=a-\left(p-1\right)^2\kappa^{2p}-p^2\kappa^{2(p-1)}=\widetilde{Q}_{p,a}(\kappa)\,.$$
As in the previous case, the monotonicity of $\widetilde{Q}_{p,a}(\kappa)$ shows that $\widetilde{Q}_{p,a}(\kappa)=0$ cannot have two positive solutions.

This finishes the proof.\hfill$\square$

\begin{rem} The result of Theorem \ref{restriction} is quite surprising when compared with the case of curves immersed in the 3-dimensional sphere $\mathbb{S}^3$. In this case, beside the critical curves for $\mathbf{\Theta}_p$ with non-constant periodic curvature, there also exist closed examples for several values of $p$ other than $p=2$ and $p\in(0,1)$ \cite{S3}.
\end{rem}

In view of above theorem and the fact that elastic curves ($p=2$) are well known \cite{LS}, it will be assumed from now on that $p\in(0,1)$. For these values, the ordinary differential equation \eqref{ODE} can be rewritten as
$$\left(\kappa'\right)^2=\frac{\kappa^2}{p^2(1-p)^2}\left(a\,\kappa^{2(1-p)}-\left[1-p\right]^2\kappa^2-p^2\right).$$
From the above expression, it is clear that if a non-constant function $\kappa(s)$ is a solution, then the right hand side must be positive, i.e.,
\begin{equation}\label{Q}
Q_{p,a}(\kappa):=a\,\kappa^{2(1-p)}-\left(1-p\right)^2\kappa^2-p^2> 0\,.
\end{equation}
This condition imposes a restriction on the parameter $a>0$. Indeed, it follows from straightforward computations that the function $Q_{p,a}(\kappa)$ has a local maximum at
$$\kappa_*:=\left(\frac{a}{1-p}\right)^{\frac{1}{2p}}\,,$$
and so $Q_{p,a}(\kappa_*)>0$ must hold. Further manipulations show that this happens precisely when 
\begin{equation}\label{d*}
a>a_*:=p^p\left(1-p\right)^{1-p}\,.
\end{equation}
Observe that when $p\to 0$ or $p\to 1$, the value of $a_*$ tends to $1$. Therefore, one can conclude that solutions to \eqref{ODE} with $p\in(0,1)$ fixed depend on two real parameters, namely $a>a_*$ and a constant of integration arising from the integration of \eqref{ODE}. Nevertheless, since the ordinary differential equation \eqref{ODE} is separable, this second constant of integration simply shifts the origin of the arc length parameter. Hence, this constant may be assumed zero after a potential translation.  Consequently, there is a (real) one-parameter family of functions $\kappa_{a}(s)$ solving \eqref{ODE} which, from the Fundamental Theorem of Curves, uniquely determine spherical curves $\gamma_{a}$ up to isometries of $\mathbb{S}^2$. Therefore, there exists a (real) one-parameter family of $p$-elastic curves $\gamma_{a}$, indexed by the real numbers $a>a_*$.

For $p\in(0,1)$, it turns out that all the solutions $\kappa_{a}(s)$ with $a>a_*$ of \eqref{ODE} are periodic functions, provided that they are defined on their maximal domain.

\begin{prop}\label{periodic} Let $p\in(0,1)$ and $a>a_*=p^p\left(1-p\right)^{1-p}$. Assume that $\gamma_{a}:I\subseteq\mathbb{R}\longrightarrow\mathbb{S}^2$ is a $p$-elastic curve with non-constant curvature defined on its maximal domain. Then, $\gamma_{a}$ is complete $(I=\mathbb{R})$ and its curvature is a periodic function.
\end{prop}
\textit{Proof.} Let $\gamma_{a}$ be a $p$-elastic curve with non-constant curvature $\kappa_{a}(s)$. Then, the curvature $\kappa_{a}(s)$ must be a solution of \eqref{ODE} for $a>a_*$. For simplicity, we introduce the new variables $u=\kappa$ and $v=u'=\kappa'$ so that \eqref{ODE} reduces to
$$v^2=\frac{u^2}{p^2(1-p)^2}\,Q_{p,a}(u)\,.$$
This is a curve representing the orbit of the differential equation \eqref{ODE} in the phase plane. Since $a>a_*$, it follows from the argument above that $Q_{p,a}(u)$ has a local maximum at $u=u_*$ and $Q_{p,a}(u_*)>0$. Combining this with $Q_{p,a}(0)=-p^2<0$ leads to the conclusion that $Q_{p,a}(u)=0$ has exactly two positive solutions, and hence this orbit must be closed. 

Moreover, the curve $C(s)=(u(s),v(s))$ is included in the trace of the above orbit and can be thought of as a bounded integral curve of the smooth vector field
$$X(u,v)=\left(v,\frac{1}{u}\left[(2-p)v^2-\frac{1}{p}u^4+\frac{1}{1-p}u^2\right]\right),$$
which is defined on $\{(u,v)\in\mathbb{R}^2\,\lvert\,u > 0\}$. Here, the upper dot denotes the derivative with respect to $u$. This implies that $C(s)$ is itself smooth and defined on the whole $\mathbb{R}$, and so are the associated curves $\gamma_{a}$. In other words, their maximal domain is the whole real line $\mathbb{R}$. Furthermore, since the vector field $X$ has no zeros along the curve $C(s)$, we conclude from the Poincar\'e-Bendixson Theorem that $C(s)$ is a periodic curve. Therefore, the non-constant curvature of $\gamma_{a}$ is a periodic function. This finishes the proof. \hfill$\square$
\\

As a consequence of Proposition \ref{periodic}, in what follows it will be assumed that $p$-elastic curves with $p\in(0,1)$ are defined on their maximal domain, i.e., on the real line $\mathbb{R}$.

\section{Existence of Closed $p$-Elastic Curves}

In this section, the existence of non-trivial closed $p$-elastic curves for $p\in(0,1)$ will be proved. Recall that a necessary, but not sufficient, condition for a curve to be closed is that its curvature is a periodic function. Consequently, Theorem \ref{restriction} shows that only the cases $p=2$ (discussed in \cite{LS}) and $p\in(0,1)$ must be considered. Therefore, our focus is placed on the case $p\in(0,1)$. For these $p$-elastic curves, the periodicity of their curvatures is guaranteed by Proposition \ref{periodic}.

Beyond periodic curvature, a suitable closure condition must also be satisfied in order to have closed curves. This condition can be explicitly obtained with the following convenient parameterization of $p$-elastic curves. Adapting the computations of \cite{LS} (see also \cite{AGM,PhD} for details) a $p$-elastic curve $\gamma_{a}$ in $\mathbb{S}^2$ can be parameterized in terms of its arc length parameter as
\begin{equation}\label{param}
\gamma_{a}(s)=\left(x_{a}(s),\sqrt{1-x_{a}^2(s)}\,\sin\psi_{a}(s),\sqrt{1-x_{a}^2(s)}\,\cos\psi_{a}(s)\right),
\end{equation}
where $x_{a}(s):=p\kappa^{p-1}(s)/\sqrt{a}$ and $\psi_{a}(s)$ is the \emph{angular progression}, which is given by the function
\begin{equation}\label{psi}
\psi_{a}(s):=(1-p)\sqrt{a\,}\int\frac{\kappa^{2-p}}{a\kappa^{2(1-p)}-p^2}\,ds\,.
\end{equation}
Observe that the definition of $x_{a}(s)$ is an extension of the well known property that elastic curves ($p=2$) can be parameterized so that one of their components is a multiple of the curvature. Moreover, if equation \eqref{ODE} is used to make a change of variable in the integral \eqref{psi}, the parameterization of $\gamma_{a}$ is given locally in terms of just one quadrature.

This parameterization can be used to deduce the geometric behavior of $p$-elastic curves in $\mathbb{S}^2$. These properties were explicitly described in \cite{MOP} for infinitely many values of $p\in(0,1)$ and can be directly extended to our general case, so they are simply stated here.

\begin{prop}\label{properties} Fix $p\in(0,1)$ and let $a>a_*=p^p\left(1-p\right)^{1-p}$. Then, the following statements hold for the $p$-elastic curve $\gamma_{a}$:
\begin{enumerate}
\item The trajectory of $\gamma_{a}$ is contained between two (proper) parallels of the half-sphere $\mathbb{S}^2_+=\{(x,y,z)\in\mathbb{S}^2\,\lvert\,x>0\}$. The curve $\gamma_{a}$ never meets the equator $x=0$ nor the pole $(1,0,0)$.
\item The curve $\gamma_{a}$ meets the bounding parallels tangentially only at the maximum and minimum values of its curvature $\kappa\equiv\kappa_{a}$.
\item The trajectory of $\gamma_{a}$ winds around the pole $(1,0,0)$ without moving backwards. In other words, its angular progression is monotonic with respect to the arc length parameter.
\item The $p$-elastic curve $\gamma_{a}$ is closed if and only if the angular progression along a period $\varrho\equiv\varrho(a)$ of the curvature is a rational multiple of $2\pi$, i.e., if and only if 
$$\Lambda_p(a):=(1-p)\sqrt{a\,}\int_0^\varrho \frac{\kappa^{2-p}}{a\kappa^{2(1-p)}-p^2}\,ds=2\pi q\,,$$
for $q\in\mathbb{Q}$.
\end{enumerate}
\end{prop}
\textit{Proof.} These statements are a straightforward extension of the properties described in \cite{MOP} to arbitrary values of $p\in(0,1)$. Similar properties can also be found in \cite{PhD} for more general functionals depending on the curvature. \hfill$\square$
\\

\begin{figure}[h!]
	\centering
	\includegraphics[height=5.55cm]{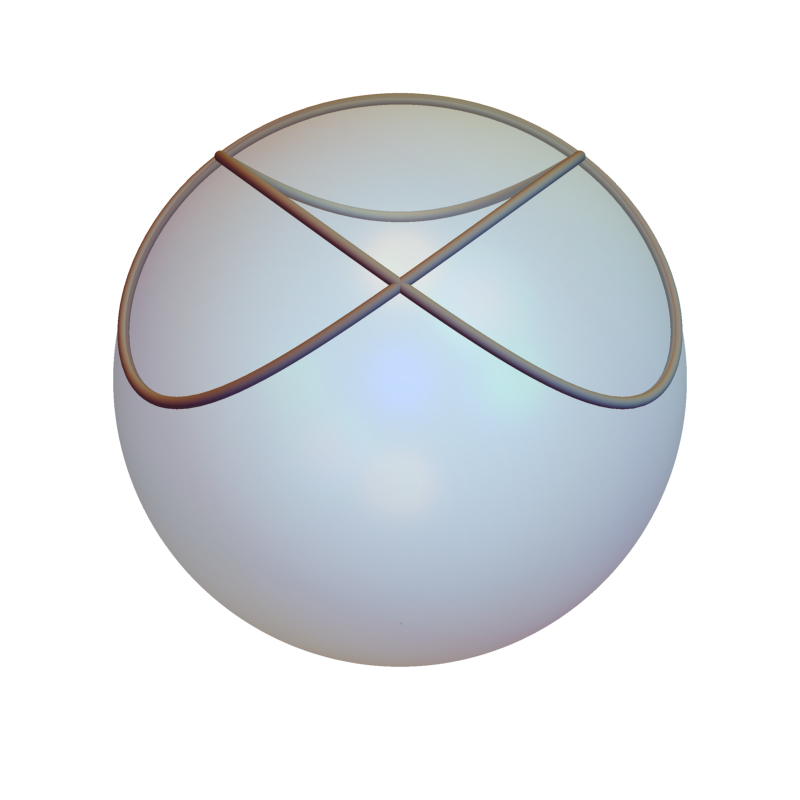}
	\includegraphics[height=5.55cm]{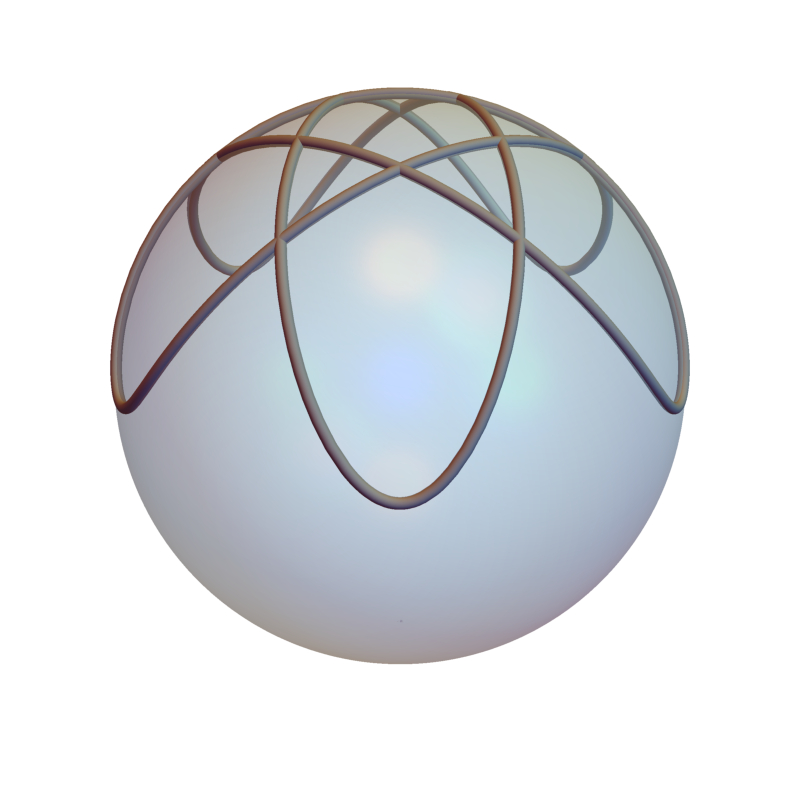}
	\includegraphics[height=5.55cm]{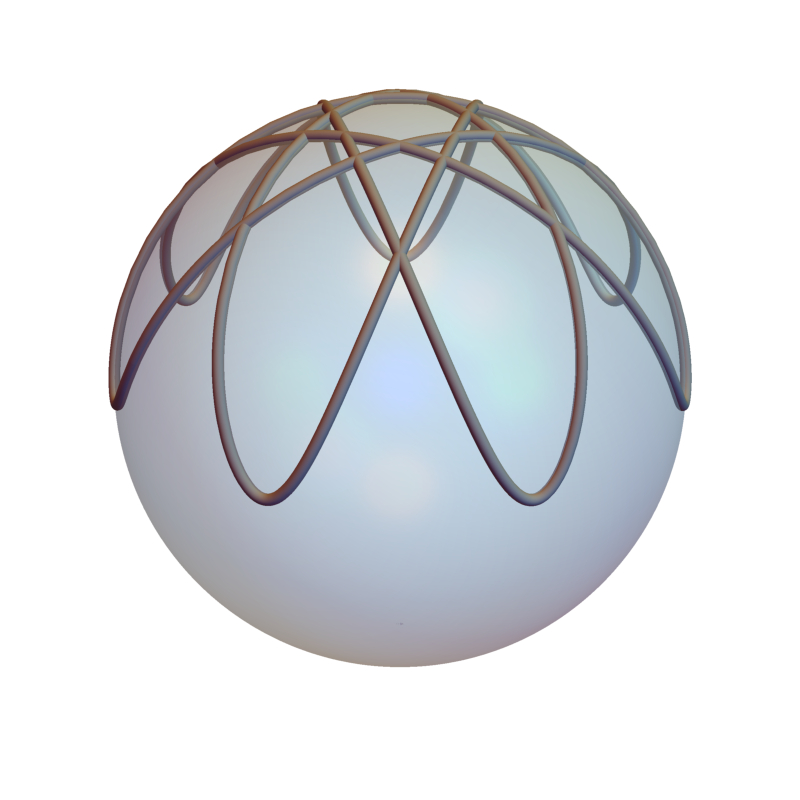}
	\caption{Three $p$-elastic curves for $p=0.3$ in $\mathbb{S}^2$ corresponding to the values $q=2/3$, $q=3/5$ and $q=4/7$, respectively. The vertical direction in the figures corresponds to the direction of the $x$-axis of $\mathbb{R}^3$.}
	\label{F1}
\end{figure}

Using the parameterization \eqref{param}, we show in Figure \ref{F1} three $p$-elastic curves for $p=0.3$ in $\mathbb{S}^2$. These curves are closed and correspond with the value of $a>a_*$ such that $\Lambda_p(a)$ is the rational multiple of $2\pi$ given by $q=2/3$, $q=3/5$ and $q=4/7$, respectively (see below for restrictions on the admissible values of $q$). In the same figure, the properties stated in Proposition \ref{properties} can be observed.

The last part of Proposition \ref{properties} suggests that $\Lambda_p(a)$ should be analyzed as a function of $a>a_*$ in order to check whether or not closed $p$-elastic curves exist. Let $\gamma_{a}$ be a $p$-elastic curve with $a>a_*$ and fixed $p\in(0,1)$. Then, from Proposition \ref{periodic}, the curvature of $\gamma_{a}$ is a non-constant periodic function. We denote by $\alpha\equiv\alpha(a)$ (respectively, $\beta\equiv\beta(a)$) the maximum (respectively, the minimum) value of the curvature of $\gamma_{a}$. Clearly, the curvature increases from $\beta$ to $\alpha$, so that in this interval \eqref{ODE} can be used to make a change of variable in $\Lambda_p(a)$, obtaining
\begin{equation}\label{Lambda}
\Lambda_p(a):=2p(1-p)^2\sqrt{a\,}\int_\beta^\alpha\frac{\kappa^{1-p}}{\left(a\kappa^{2(1-p)}-p^2\right)\sqrt{Q_{p,a}(\kappa)}}\,d\kappa\,,
\end{equation}
where $Q_{p,a}(\kappa)$ is the function defined in \eqref{Q}. The expression \eqref{Lambda} can be understood as a function $\Lambda_p:(a_*,\infty)\subset\mathbb{R}\longrightarrow\mathbb{R}$ where $a_*$ was defined in \eqref{d*}, and closed $p$-elastic curves will exist provided the image of $\Lambda_p$ is not constant. If this assertion is true, then there must exist rational multiples of $2\pi$ in its image, as desired. Moreover, it is also desirable to analyze the asymptotic behavior and the possible monotonicity of $\Lambda_p$, in order to find the possible restrictions on the rational number $q\in\mathbb{Q}$.

The asymptotic behavior of $\Lambda_p$ is analyzed in the following technical lemma.

\begin{lem}\label{lemma} Let $a_*=p^p(1-p)^{1-p}$ and $\Lambda_p:\left(a_*,\infty\right)\subset\mathbb{R}\longrightarrow\mathbb{R}$ be the function defined by the integral expression \eqref{Lambda}, for any fixed $p\in(0,1)$. Then, $\Lambda_p$ is a continuous function whose image contains the interval $(\pi,\sqrt{2}\,\pi)$.
\end{lem}
\textit{Proof.} Notice that the continuity of the function $\Lambda_p$ follows directly from the fact that $Q_{p,a}(\kappa)> 0$, since in this case,
$$a\kappa^{2(1-p)}-p^2> (1-p)^2\kappa^2>0\,.$$
Therefore, the integrand is continuous and so is $\Lambda_p$. In particular, from these inequalities and the convexity of the curves under consideration, we also deduce that the integrand in \eqref{Lambda} is positive and, hence, the image of $\Lambda_p$ is contained in the positive part of the real line.

We next study the asymptotic behavior of $\Lambda_p$. To compute the limit of $\Lambda_p(a)$ when $a\to a_*$, we apply Lemma 4.1 of \cite{Per} (see also Corollary 4.2 of the same paper). Since the function $Q_{p,a}(\kappa)$ defined in \eqref{Q} has a local maximum at $\kappa=\kappa_*$, all the conditions for this result are satisfied, and so
$$\lim_{a\to a_*}\Lambda_p(a)=2p\,(1-p)^2\sqrt{a_*}\,\frac{\kappa_*^{1-p}}{\left(a_*\kappa_*^{2(1-p)}-p^2\right)\sqrt{-\frac{1}{2}\ddot{Q}_{p,a_*}(\kappa_*)}}\,\pi\,,$$
where the upper dots denote the second derivative of $Q_{p,a_*}$ with respect to $\kappa$. A straightforward simplification using the definitions of $a_*$ and $\kappa_*$ further implies that the right hand side above is equal to $\sqrt{2}\,\pi$.

In order to compute the limit of $\Lambda_p(a)$ when $a\to \infty$, we consider first an arbitrary $p=l/t\in\mathbb{Q}\cap(0,1)$ where $l<t$ are natural numbers. We then use the change of variable $\kappa^2=u^t$ to simplify \eqref{Lambda}, obtaining
$$\Lambda_{p=l/t}(a)=l(1-p)^2\sqrt{a\,}\int_{\bar{\beta}}^{\bar{\alpha}}\frac{u^{t-1-l/2}}{\left(au^{t-l}-p^2\right)\sqrt{Q_{p,a}(u)}}\,du\,,$$
where
$$Q_{p,a}(u)=au^{t-l}-(1-p)^2u^t-p^2\,,$$
is a \emph{polynomial} of degree $t$, and $\bar{\beta}$ and $\bar{\alpha}$ are the only positive roots of $Q_{p,a}(u)$. Therefore, arguing as in Lemma 4.3 of \cite{MOP} with Cauchy's Integral Formula, we conclude that
$$\lim_{a\to\infty}\Lambda_{p}(a)=\pi\,,$$
for every $p\in\mathbb{Q}\cap(0,1)$. Finally, by the continuity of $\Lambda_p(a)$ in $p$ this result extends to every $p\in(0,1)$.

The continuity of $\Lambda_p$ in combination with the above computed limits show that the image of $\Lambda_p$ contains (not necessarily strictly) the interval $(\pi,\sqrt{2}\,\pi)$. \hfill$\square$

\begin{rem}\label{ansatz} Numerical experiments suggest that the function $\Lambda_p$ decreases from $\sqrt{2}\,\pi$ to $\pi$, and so its image is precisely the interval $(\pi,\sqrt{2}\,\pi)$. A proof of this assertion will necessitate differentiating the integral expression \eqref{Lambda} and the manipulation of complicated terms involving improper integrals (which can be seen as extensions of hyperelliptic integrals). For the case $p=1/2$ this monotonicity statement was proven in \cite{AGP1} by writing \eqref{Lambda} as a combination of standard complete elliptic integrals.
\end{rem}

From the behavior of the function $\Lambda_p$ obtained in Lemma \ref{lemma}, we deduce that for any rational number $q\in\mathbb{Q}$ satisfying $1<2q<\sqrt{2}$ there exists a value $a_q>a_*$ such that $\Lambda_p(a_q)=2\pi q$ and so, for any fixed $p\in(0,1)$, the curve $\gamma_{a_q}\equiv\gamma_q$ is a closed $p$-elastic curve. The validity of the \emph{ansatz} in Remark \ref{ansatz} would imply that for each admissible $q$ there is a curve $\gamma_q$ which is unique up to isometries of $\mathbb{S}^2$. Write the rational number $q\in\mathbb{Q}$ as $q=n/m$ for $n$ and $m$ relatively prime natural numbers. Then, the closed $p$-elastic curve $\gamma_{q}$ (also denoted $\gamma_{n,m}$) closes up in $m$ periods of its curvature, and so the trajectory of $\gamma_{n,m}$ posseses $m$ lobes and winds $n$ times around the pole $(1,0,0)$. In other words, $n$ is (perhaps up to sign) the winding number around the origin of the projection of $\gamma_{n,m}$ onto the plane $x=0$. Proving the statement of Remark \ref{ansatz} would then enable the conclusion that, for fixed $p\in(0,1)$, the family of closed $p$-elastic curves with non-constant curvature is indexed by the natural numbers $n$ and $m$ satisfying $m<2n<\sqrt{2}\,m$. In particular, (non-trivial) simple and closed $p$-elastic curves for $p\in(0,1)$ could not exist.

We summarize these findings in the following result.

\begin{thm}\label{existence} Let $n$ and $m$ be two relatively prime natural numbers satisfying $m<2n<\sqrt{2}\,m$. Then, for every $p\in(0,1)$, there exists a closed $p$-elastic curve with non-constant curvature $\gamma_{n,m}$.
\end{thm}
\textit{Proof.} Let $n$ and $m$ be two natural numbers such that ${\rm gcd}(n,m)=1$ and $m<2n<\sqrt{2}\,m$. Then, it follows that
$$\pi<2\pi\frac{n}{m}<\sqrt{2}\,\pi\,.$$
From Lemma \ref{lemma}, the function $\Lambda_p$ varies continuously from $\pi$ to $\sqrt{2}\,\pi$ and, consequently, there exists a value $a\equiv a_{n,m}>a_*$ such that
$$\Lambda_p(a_{n,m})=2\pi\frac{n}{m}\,.$$
Therefore, the last condition of Proposition \ref{properties} is satisfied for the rational number $q=n/m$ and the associated $p$-elastic curve is closed. \hfill$\square$
\\

\begin{figure}[h!]
	\centering
	\includegraphics[height=5.55cm]{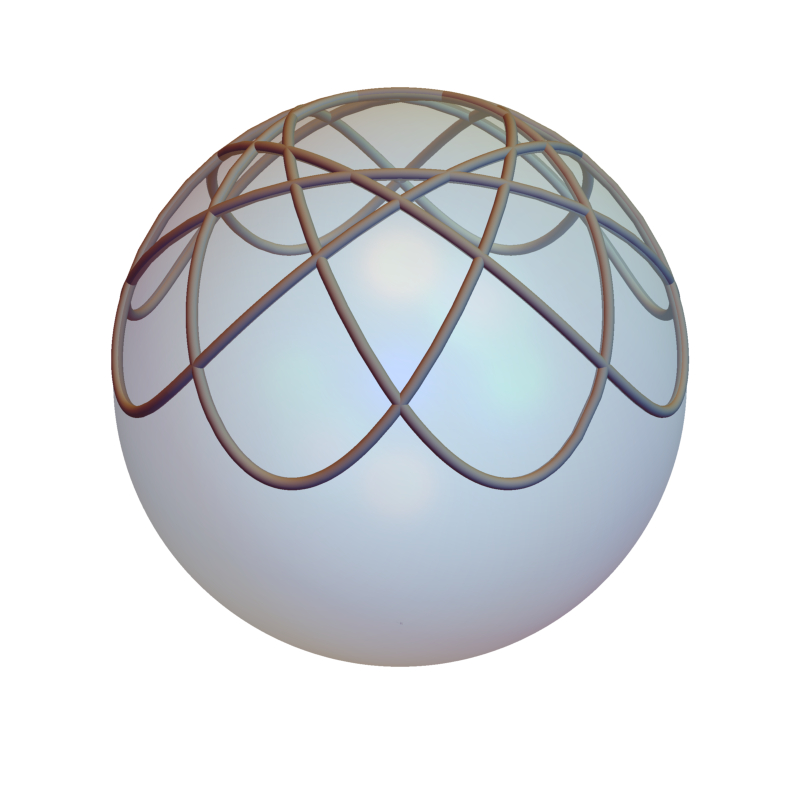}
	\includegraphics[height=5.55cm]{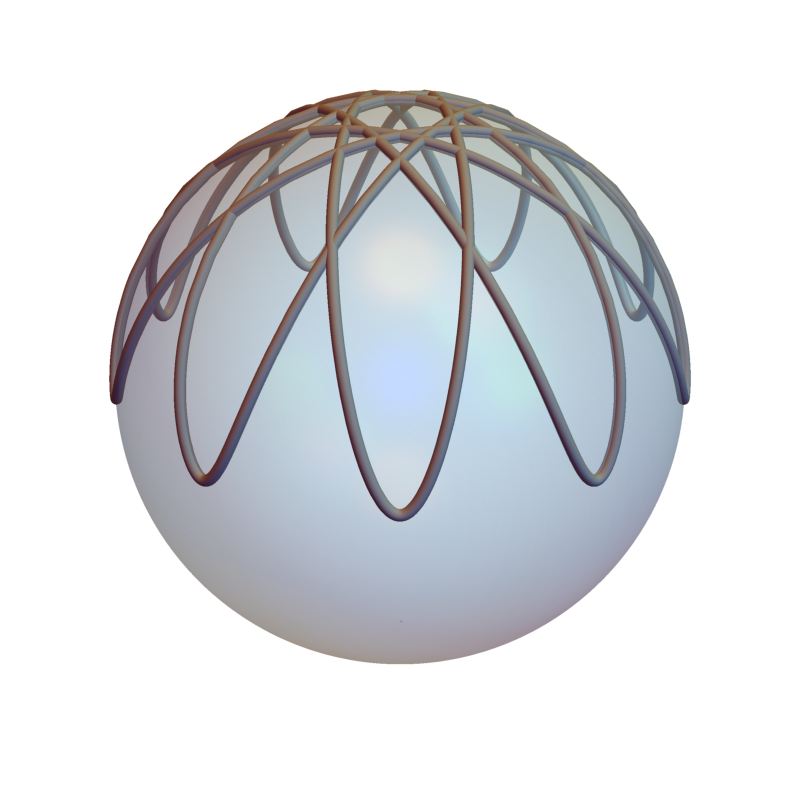}
	\includegraphics[height=5.55cm]{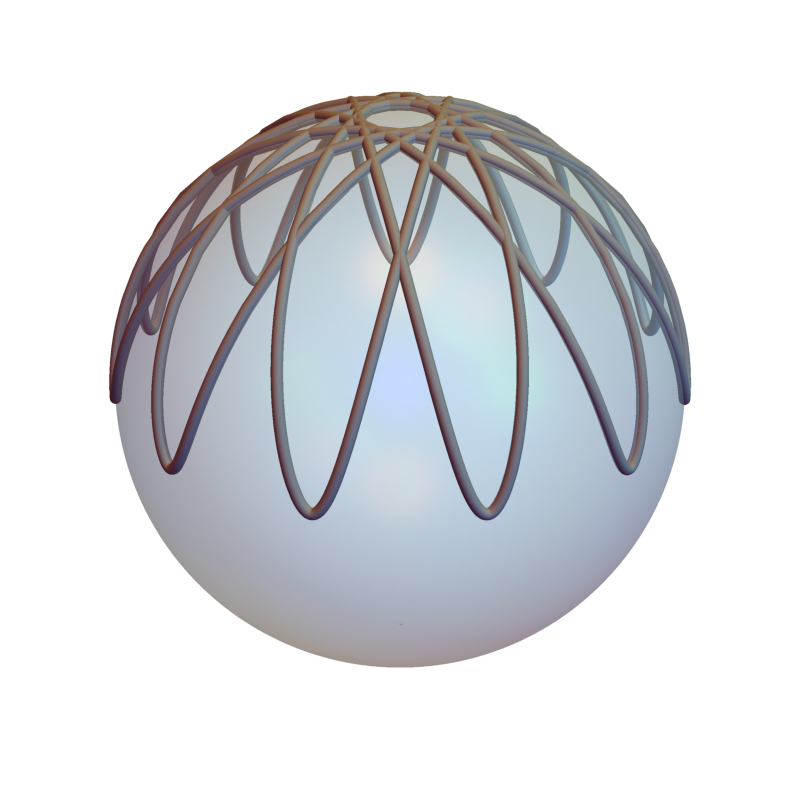}
	\caption{Closed $p$-elastic curves for $p=0.3$ in $\mathbb{S}^2$ of type $\gamma_{5,8}$, $\gamma_{5,9}$ and $\gamma_{6,11}$, respectively.}
	\label{F2}
\end{figure}

In Figure \ref{F1} we showed the ``simplest" types of closed $p$-elastic curves with non-constant curvature, for $p=0.3$. From the relation $m<2n<\sqrt{2}\,m$, these ``simplest" curves corresponded to $\gamma_{2,3}$, $\gamma_{3,5}$ and $\gamma_{4,7}$, respectively. In Figure \ref{F2} we illustrate other closed $p$-elastic curves with non-constant curvature for the same value $p=0.3$, which will correspond to $\gamma_{5,8}$, $\gamma_{5,9}$ and $\gamma_{6,11}$, respectively. Several other pictures of closed $p$-elastic curves for different values of $p$ can be found in \cite{AGM,MOP,MP,P}.

We finish this section by highlighting an interesting characteristic of closed $p$-elastic curves when $p$ varies in $(0,1)$. To build understanding, we will consider first the case of $p$-elastic circles. Recall that the Euclidean radii of $p$-elastic circles have been described in \eqref{radius}. Moreover, notice that when $p=0$  this radius is $r=1$, i.e., precisely the radius of the sphere $\mathbb{S}^2$, and on the other hand the radius $r$ approaches zero when $p\to 1$. The Euclidean radii \eqref{radius} of critical circles for $\mathbf{\Theta}_p$ with $p\in[0,1)$ can be understood as a function depending on $p$, i.e., $r\equiv r(p)$. From basic computations one can check that $r(p)$ decreases monotonically as $p$ varies from $0$ to $1$. This means that $p$-elastic circles are, up to isometries of $\mathbb{S}^2$, parallel circles which vary from the equator $x=0$ to the pole $(1,0,0)$, as $p$ increases from $0$ to $1$.

\begin{figure}[h!]
	\centering
	\includegraphics[height=3.3cm]{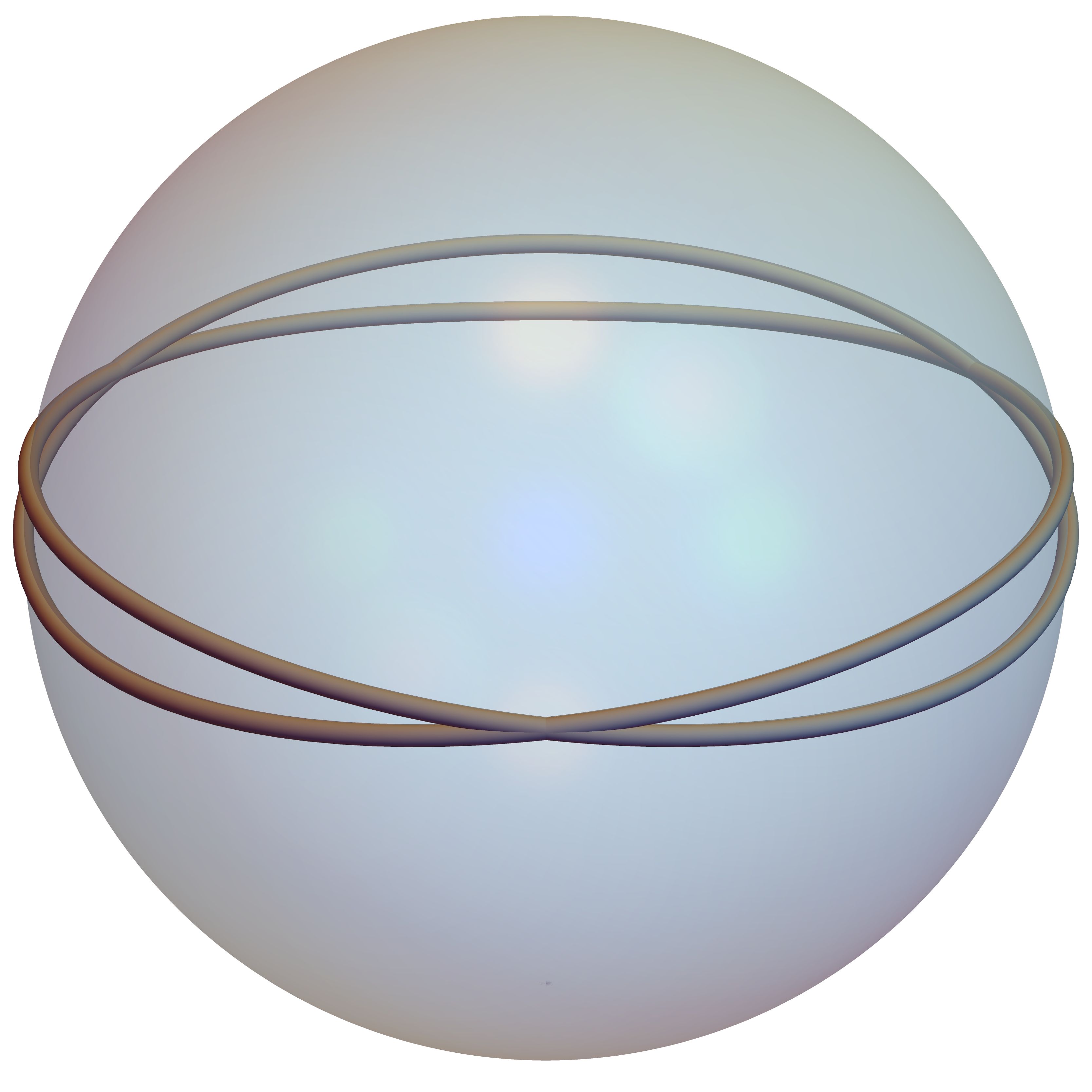}
	\includegraphics[height=3.3cm]{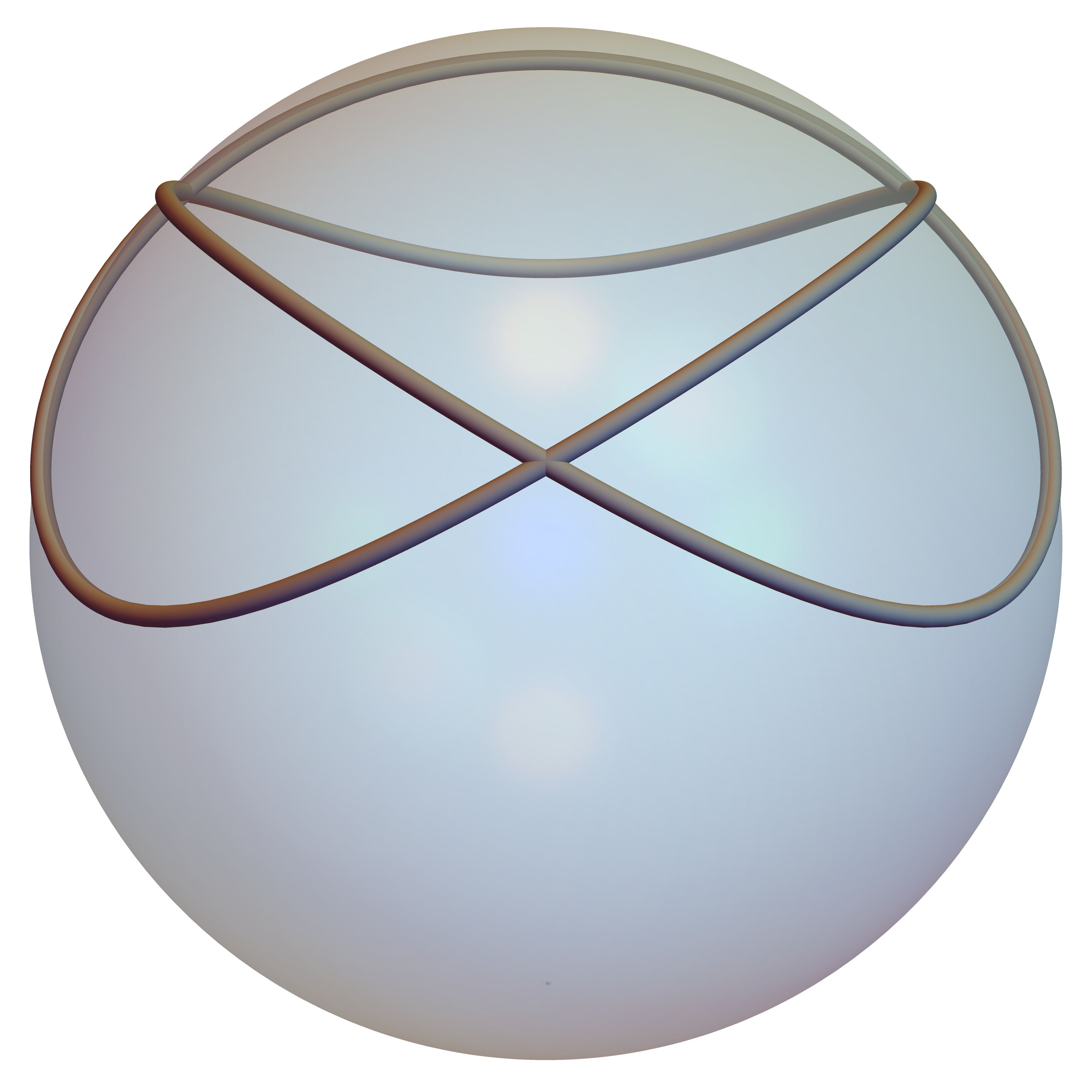}
	\includegraphics[height=3.3cm]{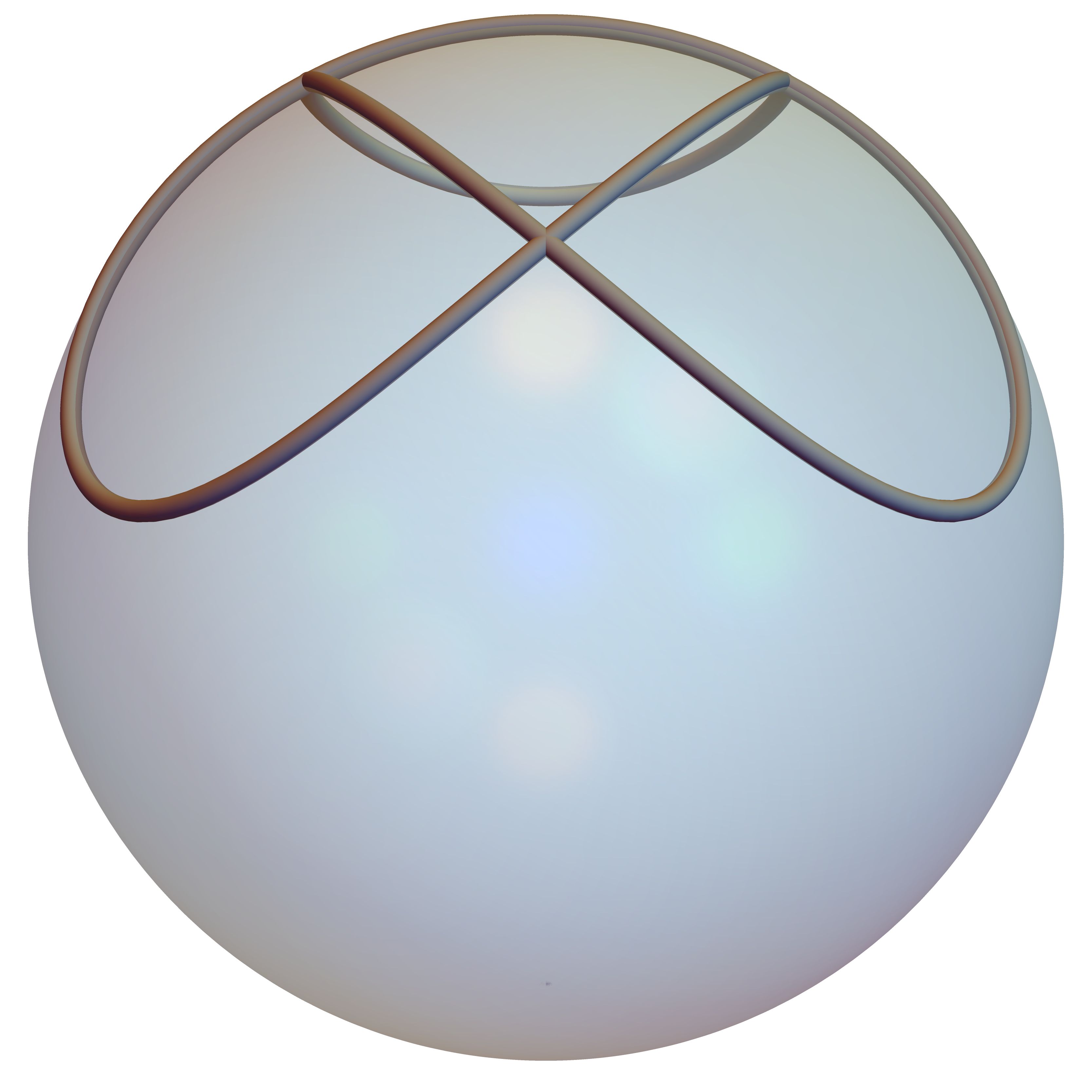}
	\includegraphics[height=3.3cm]{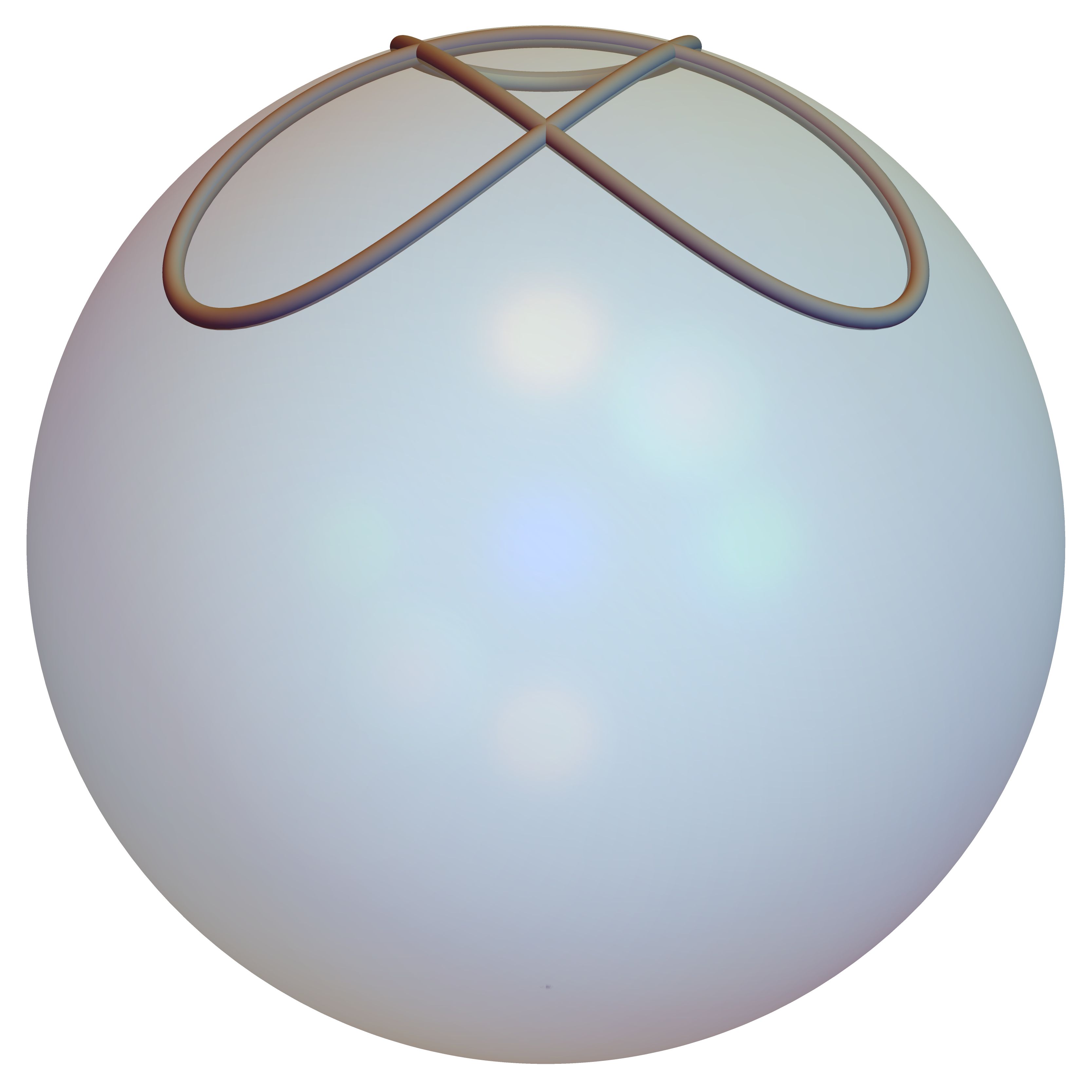}
	\includegraphics[height=3.3cm]{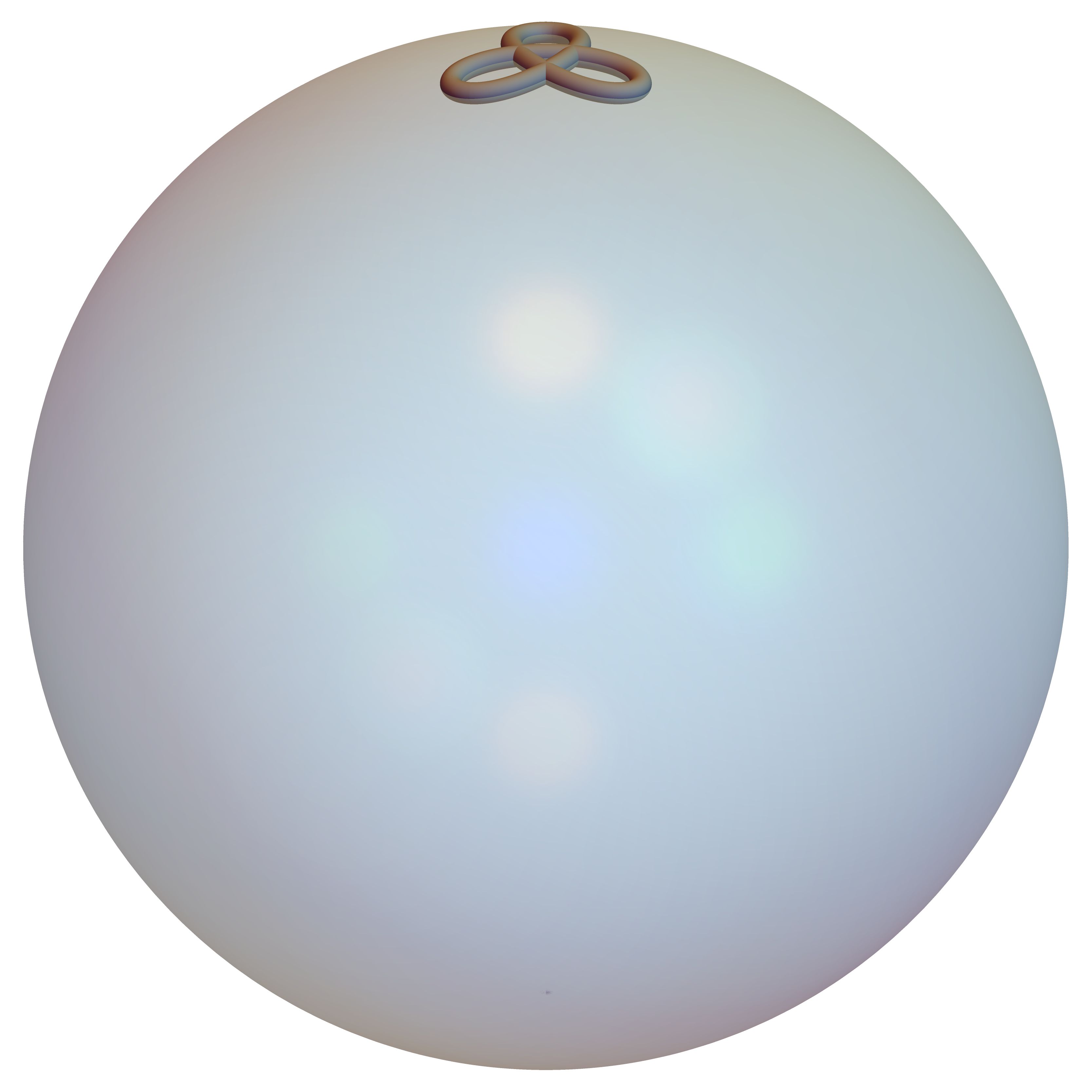}
	\caption{Closed $p$-elastic curves in $\mathbb{S}^2$ of type $\gamma_{2,3}$. From left to right: $p=0.01$, $p=0.2$, $p=0.5$, $p=0.8$ and $p=0.99$.}
	\label{motion}
\end{figure}

This phenomenon is imitated by the closed $p$-elastic curves with non-constant curvature $\gamma_{n,m}$. Indeed, their motion for fixed parameters $n$ and $m$ as $p$ varies is similar to that of the $p$-elastic circles. More precisely, let $n$ and $m$ be fixed relatively prime natural numbers satisfying $m<2n<\sqrt{2}\,m$ and let $\{\gamma_{n,m}\}_{p\in(0,1)}$ be the family of closed $p$-elastic curves $\gamma_{n,m}$ ($n$ and $m$ fixed) for different values of $p\in(0,1)$. As $p$ approaches $0$, the curves $\gamma_{n,m}$ in this family tend to the equator $x=0$, while when $p\to 1$, $\gamma_{n,m}$ shrink to the pole $(1,0,0)$. In Figure \ref{motion} we illustrate this evolution for the ``simplest'' type of closed $p$-elastic curves, i.e., for $n=2$ and $m=3$.

\section{Stability Analysis of Closed $p$-Elastic Curves}

Critical curves for $\mathbf{\Theta}_p$ are characterized by the vanishing of the first variation formula among all admissible variations. In particular, when $\mathbf{\Theta}_p$, $p\in(0,1)$, acts on the subspace of $\mathcal{C}_*^\infty(\mathbb{R},\mathbb{S}^2)$ formed by closed curves (i.e., the space of convex spherical curves which are closed, smooth, and immersed) we have shown the existence of a discrete bi-parametric family of critical curves. A natural question that arises at this point is whether or not any of these closed $p$-elastic curves is actually the minimizer of $\mathbf{\Theta}_p$, because minimizers are the only physically realizable critical points. 

We will first show that for every fixed $p\in(0,1)$ the finite infimum of $\mathbf{\Theta}_p$ cannot be attained, and so absolute minimizers of this energy functional do not exist.

\begin{thm}\label{infimum} Let $p\in(0,1)$ be fixed and assume that $\mathbf{\Theta}_p$ acts on the subspace of $\mathcal{C}^\infty_*(\mathbb{R},\mathbb{S}^2)$ formed by closed curves. Then, the infimum of $\mathbf{\Theta}_p$ is zero and is approached, but not attained, by a limit of circles.
\end{thm}
\textit{Proof.} By the definition of $\mathcal{C}_*^\infty(\mathbb{R},\mathbb{S}^2)$, our admissible curves are convex, i.e., $\kappa>0$ holds for each one. Thus, for every fixed $p\in(0,1)$, $\kappa^p>0$. It is then clear that $\mathbf{\Theta}_p(\gamma)>0$ for all closed curves in $\mathcal{C}_*^\infty(\mathbb{R},\mathbb{S}^2)$. This shows that zero is a lower bound for the energy functional $\mathbf{\Theta}_p$, which cannot be attained. We will next prove that zero is indeed the infimum.

Consider a sequence of parallels $\gamma_\epsilon$ in $\mathbb{S}^2$ of Euclidean radii $1-\epsilon$, for every $\epsilon>0$ sufficiently small. The value of $\mathbf{\Theta}_p$ at these circles satisfies the limit
$$\mathbf{\Theta}_p(\gamma_\epsilon)=\int_{\gamma_\epsilon}\kappa^p\,ds=2\pi\left(1-\epsilon\right)^{1-p}\left(1-\left[1-\epsilon\right]^2\right)^{p/2}\longrightarrow 0\,,$$
when $\epsilon\to 0$. Here we have used the relation between the curvature in $\mathbb{S}^2$ and the Euclidean radius of the circle. Therefore, zero is the infimum for $\mathbf{\Theta}_p$ and is approached by a sequence of parallels tending to an equator. \hfill$\square$

\begin{rem} Geodesics, i.e., spherical curves with $\kappa=0$ everywhere, would be absolute minimizers of $\mathbf{\Theta}_p$, $p\in(0,1)$, when acting on the space of curves with $\mathcal{L}^1$ integrable $\kappa^p$. However, they do not appear as solutions of the variational problem considered in this paper since they are not convex.
\end{rem}

Although Theorem \ref{infimum} shows that absolute minimizers of $\mathbf{\Theta}_p$ do not exist when $p\in(0,1)$, a weaker condition may also be studied. This condition is the non-negativity of the second variation formula or, in other words, the local stability of the critical curves. In \cite{LS}, it was shown that the only stable and closed (free) elastic curves ($p=2$) in $\mathbb{S}^2$ are geodesics. Motivated by this result, we now study the stability of closed $p$-elastic curves for $p\in(0,1)$.

In order to analyze the stability of $p$-elastic curves we first compute the second variation formula. For closed curves, it is enough to consider normal variations of the type $\delta\gamma=\phi(s)N(s)$, where $N(s)$ denotes the unit normal to $\gamma(s)$. Then, following the computations of Section 4 of \cite{AGM} and, after a lengthy simplification involving the Euler-Lagrange equation \eqref{EL}, the second variation formula for a $p$-elastic curve $\gamma$ reads
\begin{eqnarray*}
\delta^2\mathbf{\Theta}_p(\gamma)[\phi]&=&-p(1-p)\int_\gamma \kappa^{p-2}\left(\phi''\right)^2\,ds+(1-p)\int_\gamma \left(\left[2p+1\right]\kappa^2+2p\right)\kappa^{p-2}\left(\phi'\right)^2\,ds\\&&+\int_\gamma \mu_p(s)\phi^2\,ds\,,
\end{eqnarray*}
where
\begin{equation*}
\mu_p(s):=-p(1-p)\left(\left[p+1\right]\kappa^2+2-p\right)\kappa^{p-4}\left(\kappa'\right)^2+(1-p)\kappa^{p+2}-3\kappa^p+p\kappa^{p-2}\,.
\end{equation*}
We say that a $p$-elastic curve $\gamma$ is \emph{stable} if for every variation of $\gamma$, $\delta^2\mathbf{\Theta}_p(\gamma)\geq 0$ holds. On the contrary, we will say that the $p$-elastic curve $\gamma$ is \emph{unstable}.

As a first consequence of this second variation formula, we will conclude that $p$-elastic circles are unstable for every $p\in(0,1)$.

\begin{prop}\label{43} Let $\gamma$ be a $p$-elastic circle for $p\in(0,1)$. Then, $\gamma$ is unstable. Indeed, $\gamma$ is the absolute maximizer among once-covered circles.
\end{prop}
\textit{Proof.} To prove the instability, it is enough to find a variation of $\gamma$ such that the second variation formula in that case is negative. Here, we consider a variation of the type $\delta\gamma=N(s)$, i.e., we choose $\phi(s)\equiv 1$. Note that this variation is nothing but a perturbation of the $p$-elastic curve, which is a parallel circle, through parallels.

Since $\phi(s)\equiv 1$, both $\phi'$ and $\phi''$ are identically zero. Moreover, since the $p$-elastic curve is a circle, its curvature is constant and so $\kappa'$ is also zero. Therefore, the second variation formula simplifies to
$$\delta^2\mathbf{\Theta}_p(\gamma)[1]=\int_\gamma\left(\left[1-p\right]\kappa^4-3\kappa^2+p\right)\kappa^{p-2}\,ds=-2\frac{p}{1-p}\int_\gamma\kappa^{p-2}\,ds=-2\int_\gamma \kappa^p\,ds<0\,,$$
where, in the last two equalities, we have used that the curvature of a $p$-elastic circle is given in \eqref{circles}. This shows that $p$-elastic circles are unstable.

For the second statement, assume that $\mathbf{\Theta}_p$ is acting on the space of once-covered circles. Then, we can understand $\mathbf{\Theta}_p$ as a function of the Euclidean radius $r\in[0,1]$, which can be described as
$$\mathbf{\Theta}_p(r)=\int_\gamma\kappa^p\,ds=2\pi r\kappa^p=2\pi r^{1-p}\left(1-r^2\right)^{p/2}\,,$$
since the (geodesic) curvature $\kappa$ of a circle in $\mathbb{S}^2$ is related to its Euclidean radius $r$ by $\kappa^2=(1-r^2)/r^2$. We then use standard arguments of functions of one variable to check that the global maximum of $\mathbf{\Theta}_p(r)$ is attained at the critical point $r=\sqrt{1-p}$. \hfill$\square$

\begin{rem} Observe that the instability result also holds applying Proposition 6 of \cite{AGM}. Indeed, any possible cover of a $p$-elastic circle is unstable.
\end{rem}

We focus now on the family of closed $p$-elastic curves with non-constant curvature. For simplicity, we are going to consider a variation $\delta\gamma=N(s)$, i.e., $\phi(s)\equiv 1$. For this variation, the second variation formula simplifies to just the integral of $\mu_p(s)$ since both $\phi'$ and $\phi''$ are identically zero.

Let $\gamma\equiv\gamma_a$ be a $p$-elastic curve with non-constant curvature. Then we can use the first integral of the Euler-Lagrange equation \eqref{ODE} to get rid of $\kappa'$ in $\mu_p(s)$. Consequently, we conclude after a lengthy but  straightforward computation that, along a $p$-elastic curve with non-constant curvature and for the variation defined by $\phi(s)\equiv 1$, the second variation formula reads
\begin{eqnarray}\label{svf1}
\delta^2\mathbf{\Theta}_p(\gamma)[1]&=&\int_\gamma\left[-\frac{p+1}{p(1-p)}\,a\,\kappa^{2-p}-\frac{2-p}{p(1-p)}\,a\,\kappa^{-p}+\frac{(1-p)(2p+1)}{p}\,\kappa^{p+2}\right.\nonumber\\
&&\left.\,\,\,\quad+\frac{2(4p^2-4p+1)}{p(1-p)}\,\kappa^p+\frac{p(3-2p)}{1-p}\,\kappa^{p-2}\right]ds\,,
\end{eqnarray}
where $a>a_*=p^p(1-p)^{1-p}$ is the constant of integration appearing in \eqref{ODE}. From Proposition \ref{periodic}, it follows that the curvature of a $p$-elastic curve with non-constant curvature is a periodic function. Using the first integral \eqref{ODE} to make a change of variable in \eqref{svf1} (as done for $\Lambda_p$ in \eqref{Lambda}), we have that 
$$\delta^2\mathbf{\Theta}_p(\gamma)[1]=2m\Upsilon_p(a)\,,$$
where $m$ is the number of periods of the curvature of $\gamma_a$ and $\Upsilon_p:(a_*,\infty)\longrightarrow\mathbb{R}$ is the function defined by
$$\Upsilon_p(a):=\int_\beta^\alpha\frac{\eta_{p,a}(\kappa)}{\sqrt{Q_{p,a}(\kappa)}}\,d\kappa\,.$$
Here, $Q_{p,a}(\kappa)$ is the function defined in \eqref{Q}, $\alpha$ and $\beta$ are the only two positive solutions of $Q_{p,a}(\kappa)=0$, i.e., the maximum and minimum curvatures of $\gamma_a$, respectively, and
\begin{eqnarray*}
\eta_{p,a}(\kappa)&:=&-(p+1)a\kappa^{1-p}-(2-p)a\kappa^{-1-p}+(1-p)^2(2p+1)\kappa^{p+1}\\
&&+2(4p^2-4p+1)\kappa^{p-1}+p^2(3-2p)\kappa^{p-3}\,.
\end{eqnarray*}
Observe that $\eta_{p,a}(\kappa)$ becomes the integrand in \eqref{svf1} once we divide by $\kappa$ and multiply by $p(1-p)$, according to the change of variable \eqref{ODE}.

If $\Upsilon_p$ were to be negative, we could then conclude that  closed $p$-elastic curves are unstable for every $p\in(0,1)$. Moreover, since $\Upsilon_p$ represents the second variation of $\mathbf{\Theta}_p$ for $\phi(s)\equiv 1$ along any arch (i.e., half-period of the curvature), a strictly negative $\Upsilon_p$ would imply that any arch of a $p$-elastic curve with non-constant curvature is unstable.

One possible approach toward verifying that $\Upsilon_p$ is negative is through proving its monotonicity and asymptotic behavior. Our numerical experiments strongly support the \emph{ansatz} that $\Upsilon_p$ is a decreasing function of $a>a_*$. Moreover, we now obtain the limit when $a\to a_*$. To this end, we apply once again Corollary 4.2 of \cite{Per}. As in Lemma \ref{lemma}, the fact that $Q_{p,a}(\kappa)$ has a local maximum at $\kappa=\kappa_*$ is enough to satisfy all hypotheses of the result, and so it follows after several tedious simplifications that
\begin{equation}\label{limitUpsilon}
\lim_{a\to a_*}\Upsilon_p(a)=\frac{\eta_{p,a_*}(\kappa_*)}{\sqrt{-\frac{1}{2}\ddot{Q}_{p,a_*}(\kappa_*)}}\,\pi=-\sqrt{2a_*}\,\pi=-\sqrt{2}\,p^{p/2}(1-p)^{(1-p)/2}\pi<0\,.
\end{equation}
In conclusion, showing the \emph{ansatz} about the monotonicity of $\Upsilon_p$ will lead to a complete proof of the instability of closed $p$-elastic curves with non-constant curvature for every $p\in(0,1)$. However, the proof of the \emph{ansatz} would involve similar computations to those explained in Remark \ref{ansatz}. Therefore, although the idea is straightforward, a rigorous approach in this way may seem prohibitively difficult. Nevertheless, using a trick involving integration by parts we manage to prove the instability in a different way.

\begin{thm}\label{instability} For every $p\in(0,1)$, closed $p$-elastic curves with non-constant curvature are unstable.
\end{thm}
\textit{Proof.} We begin by obtaining a relation between integrals of the type arising in the definition of $\Upsilon_p$. Since $\alpha$ and $\beta$ are the positive solutions of $Q_{p,a}(\kappa)=0$, using the classical formula of integration by parts, we have that
$$0=\left(\kappa^t\sqrt{Q_{p,a}(\kappa)}\right)_{\beta}^\alpha=\int_\beta^\alpha\frac{\frac{1}{2}\kappa^t\dot{Q}_{p,a}(\kappa)+t\kappa^{t-1}Q_{p,a}(\kappa)}{\sqrt{Q_{p,a}(\kappa)}}\,d\kappa\,,$$
for every $t\in\mathbb{R}$. Next, from the particular expression of $Q_{p,a}(\kappa)$ given in \eqref{Q} we conclude that above equality reads
\begin{equation}\label{parts}
	(1-p)^2(1+t)\int_\beta^\alpha \frac{\kappa^{1+t}}{\sqrt{Q_{p,a}(\kappa)}}\,d\kappa=a(1+t-p)\int_\beta^\alpha\frac{\kappa^{1+t-2p}}{\sqrt{Q_{p,a}(\kappa)}}\,d\kappa-tp^2\int_\beta^\alpha\frac{\kappa^{-1+t}}{\sqrt{Q_{p,a}(\kappa)}}\,d\kappa\,.
\end{equation}
We will now use this relation to simplify the expression of $\Upsilon_p$. Using it with $t=p$ in the third integral of $\Upsilon_p$ and with $t=p-2$ in the last one, we obtain after several simplifications that
\begin{eqnarray}
\Upsilon_p(a)&=&-\frac{a p^2}{1+p}\int_\beta^\alpha\frac{\kappa^{1-p}}{\sqrt{Q_{p,a}(\kappa)}}\,d\kappa-\frac{a(1-p)^2}{2-p}\int_\beta^\alpha\frac{\kappa^{-1-p}}{\sqrt{Q_{p,a}(\kappa)}}\,d\kappa\nonumber\\
&&+\frac{-4p^4+8p^3+2p^2-6p+1}{(1+p)(2-p)}\int_\beta^\alpha\frac{\kappa^{-1+p}}{\sqrt{Q_{p,a}(\kappa)}}\,d\kappa\,.\label{UPSILON}
\end{eqnarray} 
Since the curves under consideration are all convex, it follows that all the integrals in above expression are positive. Furthermore, the coefficients of the first two integrals are negative for every $p\in(0,1)$, while the coefficient of the third one is non-positive provided that $p\in[(1-\sqrt{4-\sqrt{13}})/2,(1+\sqrt{4-\sqrt{13}})/2]$. Therefore, for these values of $p$, $\Upsilon_p(a)<0$ and so $\delta^2\mathbf{\Theta}_p(\gamma)[1]=2m\Upsilon_p(a)<0$, yielding the conclusion. 

We assume now that $p\in(1+\sqrt{4-\sqrt{13}})/2,1)$. In this case, $-4p^4+8p^3+2p^2-6p+1>0$ holds, which means that the third term in \eqref{UPSILON} is positive. We then use \eqref{parts} for $t=p$ to rewrite this last integral, obtaining a different expression for $\Upsilon_p$, namely,
\begin{eqnarray*}
	\Upsilon_p(a)&=&-\frac{a(1-p)(p^4-2p^3-3p^2+6p-1)}{p^3(2-p)}\int_\beta^\alpha\frac{\kappa^{1-p}}{\sqrt{Q_{p,a}(\kappa)}}\,d\kappa-\frac{a(1-p)^2}{2-p}\int_\beta^\alpha\frac{\kappa^{-1-p}}{\sqrt{Q_{p,a}(\kappa)}}\,d\kappa\\
	&&-\frac{(1-p)^2(-4p^4+8p^3+2p^2-6p+1)}{p^3(2-p)}\int_\beta^\alpha\frac{\kappa^{1+p}}{\sqrt{Q_{p,a}(\kappa)}}\,d\kappa\,.
\end{eqnarray*}
As above, we deduce from the convexity of our curves that the three integrals are positive. The coefficient of the second integral is clearly negative for every $p\in(0,1)$ and so is the one of the third term, because the inequality $-4p^4+8p^3+2p^2-6p+1>0$ holds. A simple analysis also shows that $p^4-2p^3-3p^2+6p-1>0$ when $p\in(1+\sqrt{4-\sqrt{13}})/2,1)$. We then conclude that, in this case, $\Upsilon_p(a)<0$ and so $\delta^2\mathbf{\Theta}_p(\gamma)[1]=2m\Upsilon_p(a)<0$.

Finally, take $p\in(0,(1-\sqrt{4-\sqrt{13}})/2)$. As before, $-4p^4+8p^3+2p^2-6p+1>0$ holds and so we need to rewrite the last integral in \eqref{UPSILON}. Applying again \eqref{parts} with $t=p-2$, we derive an expression for $\Upsilon_p$ as
\begin{eqnarray*}
	\Upsilon_p(a)&=&-\frac{a p^2}{1+p}\int_\beta^\alpha\frac{\kappa^{1-p}}{\sqrt{Q_{p,a}(\kappa)}}\,d\kappa-\frac{p^2(-4p^4+8p^3+2p^2-6p+1)}{(1+p)(1-p)^3}\int_\beta^\alpha\frac{\kappa^{-3+p}}{\sqrt{Q_{p,a}(\kappa)}}\,d\kappa\\
	&&-\frac{a p(-p^5+4p^4-p^3-8p^2+3p+2)}{(1+p)(1-p)^3(2-p)}\int_\beta^\alpha\frac{\kappa^{-1-p}}{\sqrt{Q_{p,a}(\kappa)}}\,d\kappa\,.
\end{eqnarray*}
As in previous cases, the integrals are positive while the coefficients are negative, and so $\Upsilon_p(a)<0$ for these values of $p$, too.

Consequently, for every $p\in(0,1)$, closed $p$-elastic curves with non-constant curvature are unstable. \hfill$\square$

\begin{rem} In fact, the proof of Theorem \ref{instability} shows a stronger result. Indeed, for $p\in(0,1)$, any arch of a $p$-elastic curve with non-constant curvature is unstable.
\end{rem}

In the particular case $p=1/2$, the function $Q_{1/2,a}$ is a second order \emph{polynomial}, hence the function $\Upsilon_{1/2}(a)$ can be written as a combination of standard complete elliptic integrals of the first and second kind, denoted by $K(\zeta)$ and $E(\zeta)$, respectively. This simplifies the computations considerably and gives insight into the second variation formula.

\begin{thm}\label{stab12} Let $\gamma$ be a closed $1/2$-elastic curve closing up in $m$ periods of the curvature. Then,
$$\delta^2\mathbf{\Theta}_{1/2}(\gamma)[1]=2m\Upsilon_{1/2}(a)=-\frac{8}{3}m\left(\sqrt{\alpha}\,a\, E(\zeta)+\sqrt{\beta}\,K(\zeta)\right)<0\,,$$
where $\zeta:=\sqrt{(\alpha-\beta)/\alpha}$ and $\alpha>\beta$ are the positive roots of $Q_{1/2,a}(\kappa)=0$. Consequently, closed $1/2$-elastic curves are unstable.
\end{thm}
\textit{Proof.} Consider $p=1/2$. Plugging this in the definitions of $\eta_{1/2,a}$, $Q_{1/2,a}$, and splitting the integral $\Upsilon_{1/2}(a)$ into four simpler integrals, we have that
\begin{eqnarray}\label{split}
\Upsilon_{1/2}(a)&=&-3a\int_\beta^\alpha \frac{\kappa}{\sqrt{(\alpha-\kappa)(\kappa-\beta)\kappa\,}}\,d\kappa-3a\int_\beta^\alpha \frac{1}{\sqrt{(\alpha-\kappa)(\kappa-\beta)\kappa^3\,}}\,d\kappa\nonumber\\
&&+\int_\beta^\alpha \frac{\kappa^2}{\sqrt{(\alpha-\kappa)(\kappa-\beta)\kappa\,}}\,d\kappa+\int_\beta^\alpha \frac{1}{\sqrt{(\alpha-\kappa)(\kappa-\beta)\kappa^5\,}}\,d\kappa\,.
\end{eqnarray}
Each of the integrals above can be rewritten in terms of complete elliptic integrals of first and second kind. Recall that these complete elliptic integrals are defined, respectively, by
$$K(\zeta):=\int_0^1 \frac{1}{\sqrt{(1-u^2)(1-\zeta^2u^2)}}\,du\,,\quad\quad\text{and}\quad\quad E(\zeta):=\int_0^1\frac{\sqrt{1-\zeta^2u^2}}{\sqrt{1-u^2}}\,du\,.$$
Using the corresponding formulas of \cite{Gr}, the first integral in \eqref{split} is $2\sqrt{\alpha}E(\zeta)$ (see 3.132-5), where $\zeta:=\sqrt{(\alpha-\beta)/\alpha}$. Similarly, from 3.133-16, the second integral in \eqref{split} is again $2\sqrt{\alpha}E(\zeta)$; while, it follows from 3.134-16 that the fourth one is
$$\int_\beta^\alpha \frac{1}{\sqrt{(\alpha-\kappa)(\kappa-\beta)\kappa^5\,}}\,d\kappa=\frac{2\sqrt{\alpha}}{3}\left(2\left[\alpha+\beta\right]E(\zeta)-\beta K(\zeta)\right).$$
For the third integral in \eqref{split} we make the change of variable $\kappa=u^2$, so that
$$\int_\beta^\alpha \frac{\kappa^2}{\sqrt{(\alpha-\kappa)(\kappa-\beta)\kappa\,}}\,d\kappa=2\int_{\sqrt{\beta}}^{\sqrt{\alpha}}\frac{u^4}{\sqrt{(\alpha-u^2)(u^2-\beta)}}\,du\,,$$
which, according to 3.154-7, simplifies to
$$\frac{4}{3}\sqrt{\alpha}\left(\alpha+\beta\right)E(\zeta)-\frac{2}{3}\sqrt{\beta}K(\zeta)\,.$$
In all these computations we have used the Cardano-Vieta relations to obtain that $\alpha\beta=1$ and $\alpha+\beta=4a$  which leads to a further simplification of the result. Combining everything, we conclude that
\begin{equation}\label{ups}
\Upsilon_{1/2}(a)=-\frac{4}{3}\sqrt{\alpha}\,a\,E(\zeta)-\frac{4}{3}\sqrt{\beta}\,K(\zeta)\,.
\end{equation} 
The first statement then follows since $\Upsilon_{1/2}(a)$ represents the second variation for $\phi(s)\equiv 1$ along one half-period of the curvature.

For the second statement, observe that $\Upsilon_{1/2}(a)$ given in \eqref{ups} is negative and so is
$$\delta^2\mathbf{\Theta}_{1/2}(\gamma)[1]=2m\Upsilon_{1/2}(a)\,.$$
Therefore, for every closed $1/2$-elastic curve with non-constant curvature there exists a variation (namely $\delta\gamma=N(s)$) such that $\delta^2\mathbf{\Theta}_{1/2}(\gamma)[\phi]<0$, and  hence $\gamma$ is unstable. \hfill$\square$

\begin{rem} In the expression of $\Upsilon_{1/2}(a)$ of Theorem \ref{stab12}, we can take the limit when $a\to a_*=1/2$. Note that both $\beta$ and $\alpha$ tend to $1$ when $a\to a_*$, and so $\zeta\to 0$. It then follows from $K(0)=E(0)=\pi/2$ that
$$\lim_{a\to a_*=1/2} \Upsilon_{1/2}(a)=-\pi\,,$$
which coincides with the general limit computed in \eqref{limitUpsilon}. Moreover, one can check by differentiating \eqref{ups} with respect to $a$ and using the classical formulas for the derivatives of the elliptic integrals that $\Upsilon_{1/2}(a)$ is decreasing, as claimed above for the general case.
\end{rem}

\section*{Appendix A. Numerical Computations}

Computations in Table \ref{tab1} show the numerical values associated to the closed $p$-elastic curves illustrated throughout the paper. For each figure we first show the value of the energy parameter $p$ and the type of the curve $\gamma_{n,m}$ (in other words, the number of winds around the pole $n$, and the number of periods of the curvature, i.e., lobes, $m$). Recall that these numbers satisfy $m<2n<\sqrt{2}\,m$.

Moreover, in the same table we also give the numerical value of $a\equiv a_{n,m}$ for the corresponding curve, that is, the numerical solution of 
$$\Lambda_p(a)=2\pi q=2\pi\frac{n}{m}\,,$$
for fixed $n$ and $m$ (see Proposition \ref{properties} and Theorem \ref{existence} for details). For the examples with fixed $p=0.3$, it can be checked that the numerical values of $a\equiv a_{n,m}$ decrease as the quotient $n/m$ increases, an expected result from our discussion of the \emph{conjectured} monotonicity of $\Lambda_p$ (see Remark \ref{ansatz}).

We also compute the value of the energy functional $\mathbf{\Theta}_p$ along each of the $p$-elastic curves. For a closed $p$-elastic curve $\gamma_{n,m}\equiv \gamma_a$ closing up in $m$ periods of its curvature, this value is computed by numerically solving the integral
$$\mathbf{\Theta}_p(\gamma_a)=\int_\gamma \kappa^p\, ds=2\, m\, p\left(1-p\right)\int_\beta^\alpha \frac{\kappa^{p-1}}{\sqrt{Q_{p,a}(\kappa)}}\,d\kappa\,,$$
where $\alpha>\beta$ are the maximum and minimum curvatures of $\gamma_{n,m}$. In other words, the only positive solutions of $Q_{p,a}(\kappa)=0$. The limit when $a\to a_*$ of above expression can be computed analytically. Using the result of \cite{Per} as done in the first part of the proof of Lemma \ref{lemma} and \eqref{limitUpsilon}, we obtain
$$\lim_{a\to a_*}\mathbf{\Theta}_p(\gamma_a)=2 \,m\, p \left(1-p\right) \frac{\kappa_*^{p-1}}{\sqrt{-\frac{1}{2}\ddot{Q}_{p,a_*}(\kappa_*)}}=m\sqrt{2a_*}\,\pi=m\sqrt{2}\,p^{p/2}(1-p)^{(1-p)/2}\pi\,.$$
Observe that, comparing this limit with the one in \eqref{limitUpsilon}, it can be concluded that $\delta^2\mathbf{\Theta}_p(\gamma)[1]=-2\mathbf{\Theta}_p(\gamma)$ for a $p$-elastic circle $\gamma$. This agrees with the computation of the proof of Proposition \ref{43}. In order to compare with our theoretical findings, we check that the numerical values of the energy obtained in Table \ref{tab1} approach this limit as the value of $a\equiv a_{n,m}$ decreases (for a fixed value of $p$). We also recall here that the infimum of $\mathbf{\Theta}_p$ with $p\in(0,1)$ is zero (Theorem \ref{infimum}). Our experimental values are far from this infimum, which suggests that the closed $p$-elastic curves are far from minimizing $\mathbf{\Theta}_p$. Indeed, we also obtain the numerical value of $\delta^2\mathbf{\Theta}_p(\gamma)[1]$, i.e., of the second variation for $\delta\gamma=N(s)$. All these numbers are negative and so the closed $p$-elastic curves are unstable, as proved in Theorem \ref{instability}. We can moreover observe that the further $a\equiv a_{n,m}$ is from $a_*$, the more negative this number becomes, which follows from the \emph{conjectured} monotonicity of $\Upsilon_p(a)$.

\begin{center}
\begin{table}[h!]
\caption{Numerical values associated to the closed $p$-elastic curves of Figures \ref{F1}-\ref{motion}. All numbers are rounded up to just two decimals.}\label{tab1}
\makebox[\textwidth][c]{
\begin{tabular}{|c|c|c|c|c|c|} 
\hline
Figure & Value of $p$ & Curve Type ($\gamma_{n,m}$) & Value of $a$ & Energy ($\mathbf{\Theta}_p$) & $\delta^2\mathbf{\Theta}_p(\gamma)[1]$ \\
\hline
Fig. \ref{F1} (Left) & $0.3$ & $\gamma_{2,3}$ & $0.79$ & $9.2$ & $-24.88$ \\
Fig. \ref{F1} (Center) & $0.3$ & $\gamma_{3,5}$ & $1.68$ & $13.1$ & $-85.05$ \\
Fig. \ref{F1} (Right) & $0.3$ & $\gamma_{4,7}$ & $2.66$ & $16.53$ & $-137.51$ \\
Fig. \ref{F2} (Left) & $0.3$ & $\gamma_{5,8}$ & $1.23$ & $22.87$ & $-92.17$ \\
Fig. \ref{F2} (Center) & $0.3$ & $\gamma_{5,9}$ & $3.74$ & $19.65$ & $-451.26$ \\
Fig. \ref{F2} (Right) & $0.3$ & $\gamma_{6,11}$ & $4.9$ & $22.57$ & $-841.27$ \\
Fig. \ref{motion} (Left) & $0.01$ & $\gamma_{2,3}$ & $0.96$ & $12.22$ & $-33.34$ \\
Fig. \ref{motion} (Center-Left) & $0.2$ & $\gamma_{2,3}$ & $0.79$ & $9.74$ & $-26.37$ \\
Fig. \ref{motion} (Center) & $0.5$ & $\gamma_{2,3}$ & $0.8$ & $8.82$ & $-23.88$ \\
Fig. \ref{motion} (Center-Right) & $0.8$ & $\gamma_{2,3}$ & $0.79$ & $9.74$ & $-26.37$ \\
Fig. \ref{motion} (Right) & $0.99$ & $\gamma_{2,3}$ & $0.96$ & $12.22$ & $-33.34$ \\
\hline
\end{tabular}}
\end{table}
\end{center}

\section*{Appendix B. Application to $p$-Willmore Tori}

In this appendix we show an application of the theory of spherical $p$-elastic curves to the existence of closed $p$-Willmore surfaces in spheres. For convenience, we consider without loss of generality the $3$-dimensional sphere of radius two 
$$\mathbb{S}^3(1/4)=\{(x,y,z,t)\in\mathbb{R}^4\,\lvert\,x^2+y^2+z^2+t^2=4\}\subset\mathbb{R}^4\,,$$
where $(x,y,z,t)$ are the standard coordinates of the Euclidean space $\mathbb{R}^4$. For simplicity we will denote $\mathbb{S}^3(1/4)$ by $\mathbb{S}^3$.

Let $\Sigma$ be an oriented, connected and closed (i.e., compact without boundary) surface and consider a smooth immersion $X:\Sigma\longrightarrow\mathbb{S}^3$. We denote by $\mathcal{C}^\infty(\Sigma,\mathbb{S}^3)$ the space of these immersions. An immersion $X\in\mathcal{C}^\infty(\Sigma,\mathbb{S}^3)$ is said to be \emph{weakly convex} if its Gaussian curvature $K$ is non-negative and its mean curvature $H$ is strictly positive with respect to the orientation of $\Sigma$. The subspace of weakly convex immersions will be denoted by $\mathcal{C}_*^\infty(\Sigma,\mathbb{S}^3)$. We say that $X:\Sigma\longrightarrow\mathbb{S}^3$ is a \emph{$p$-Willmore surface} (or, generalized Willmore surface) if it is a critical point of the functional,
\begin{equation}\label{energysurface}
	\mathcal{W}_p(X):=\int_\Sigma H^p\,dA\,,
\end{equation}
where $p\in\mathbb{R}$ and $dA$ denotes the area element. The functional $\mathcal{W}_p$ will be assumed to act on the space $\mathcal{C}^\infty(\Sigma,\mathbb{S}^3)$ if $p\in\mathbb{N}$ (recall that $0\in\mathbb{N}$), or on $\mathcal{C}^\infty_*(\Sigma,\mathbb{S}^3)$ if $p\in\mathbb{R}\setminus\mathbb{N}$. This family of functionals include some classical ones such as the area functional ($p=0$) and the Willmore energy ($p=2$). For more information about these functionals and a historical background see, for instance, \cite{Anthony} and references therein. From standard arguments of the Calculus of Variations (for details, see for instance \cite{Anthony,P}), we can obtain the Euler-Lagrange equation associated to $\mathcal{W}_p$, which is
\begin{equation}\label{ELsurface}
	p\Delta\left(H^{p-1}\right)+4\left(p-1\right)H^{p+1}+p\left(1-2K\right)H^{p-1}=0\,,
\end{equation}
where $\Delta$ is the Laplacian of the surface. Equation \eqref{ELsurface} is, in general, a fourth order non-linear elliptic partial differential equation, and so obtaining global solutions is a difficult task. Solutions of this equation have been long studied for the area functional ($p=0$), in which case \eqref{ELsurface} is of second order, and also for the Willmore energy ($p=2$) where the conformal invariance is a powerful ally. For other choices of the parameter $p$, there are not many known closed solutions. We point out here that critical Hopf tori for the case $p=1/2$ were obtained in \cite{P}.

Employing the Symmetric Criticality Principle of Palais \cite{Pa} and the Hopf tori constructed by Pinkall \cite{Pi}, we will next obtain solutions of \eqref{ELsurface}, based on $p$-elastic curves. We begin by briefly describing the construction of Hopf tori, which is slightly different from the one introduced in \cite{Pi} (for details, see \cite{PhD} and references therein). For the sake of simplicity it is convenient to identify $\mathbb{R}^2$ with the complex plane $\mathbb{C}$ and so $\mathbb{R}^4$ with $\mathbb{C}^2$. In $\mathbb{C}^2$ we are going to consider the Riemannian metric
$$g\left((z_1,z_2),(\omega_1,\omega_2)\right)=\Re\left(z_1\overline{\omega}_1+z_2\overline{\omega}_2\right),$$
where $z_i$ and $\omega_i$, $i=1,2$ are complex numbers and the bar denotes the usual complex conjugate. Consequently, in $(\mathbb{C}^2,g)$ the $3$-dimensional sphere of radius two $\mathbb{S}^3$ can be seen as the hyperquadric $g\left((z_1,z_2),(\omega_1,\omega_2)\right)=4$. Define the map $\widetilde{\pi}:\mathbb{C}^2\longrightarrow\mathbb{C}^2$ by
$$\widetilde{\pi}(z,\omega)=\frac{1}{4}\left(\lvert z\rvert^2-\lvert\omega\rvert^2,2\overline{z}\omega\right).$$
Then, the restriction of $\widetilde{\pi}$ (also denoted by $\widetilde{\pi}$) to $\mathbb{S}^3$, identified as the above described hyperquadric, gives the \emph{standard Hopf mapping} $\widetilde{\pi}:\mathbb{S}^3\equiv\mathbb{S}^3(1/4)\longrightarrow\mathbb{S}^2\equiv\mathbb{S}^2(1)$.

Let $\gamma$ be a curve in $\mathbb{S}^2$, then the complete lift $X_\gamma:=\widetilde{\pi}^{-1}(\gamma)$ is a flat surface in $\mathbb{S}^3$ whose mean curvature $H$ satisfies
\begin{equation}\label{H}
	H=\frac{1}{2}\left(\kappa\circ \widetilde{\pi}\right),
\end{equation}
where $\kappa$ denotes the geodesic curvature of $\gamma$ in $\mathbb{S}^2$. The surface $X_\gamma$ is usually called the \emph{Hopf cylinder} (or, Hopf tube) based on $\gamma$. Observe that the fibers of $X_\gamma$ are isomorphic to $\mathbb{S}^1$. Therefore, if the curve $\gamma$ is closed so is $X_\gamma$ and we refer to it as the \emph{Hopf torus} based on $\gamma$. The covering map
$$X_\gamma(t,s):=e^{it}\overline{\gamma}(s)\,,$$
where $\overline{\gamma}$ denotes a \emph{horizontal lift} of $\gamma$ can be used to parameterize the Hopf cylinder. Explicit parameterizations of horizontal lifts $\overline{\gamma}$ can be found, for instance, in \cite{P,PhD}. Note that even if the curve $\gamma$ is closed, its horizontal lift may not be closed due to the non-trivial holonomy. However, if in addition the area enclosed by $\gamma$ in $\mathbb{S}^2$ is a rational multiple of $\pi$, then the horizontal lift of a suitable cover of $\gamma$ will be closed \cite{ABG0}.

Let $\gamma$ be a closed curve in $\mathbb{S}^2$. It follows from the Symmetric Criticality Principle of Palais \cite{Pa} that the Hopf torus $X_\gamma$ is a $p$-Willmore surface in $\mathbb{S}^3$ if and only if $\gamma$ is a $p$-elastic curve. Although the proof of this statement relies heavily on the Symmetric Criticality Principle, a direct relation between the Euler-Lagrange equations \eqref{EL} and \eqref{ELsurface} can be obtained using the expression of the mean curvature of $X_\gamma$ given in \eqref{H} and the fact that $X_\gamma$ is flat, i.e., its Gaussian curvature is zero ($K=0$).

In any case, we have that $p$-Willmore Hopf tori in $\mathbb{S}^3$ exist provided that closed $p$-elastic curves in $\mathbb{S}^2$ exist. From Theorem \ref{restriction}, this occurs only in the cases $p=2$ and $p\in(0,1)$. The case $p=2$ was studied in \cite{Pi} and, in combination with \cite{LS}, the existence of infinitely many non-conformally minimal Willmore tori was shown. Here, we focus on the case $p\in(0,1)$.  Translating our findings of Theorem \ref{existence} for $p$-elastic curves, we conclude that for every $p\in(0,1)$ and any pair of relatively prime natural numbers $(n,m)$ satisfying $m<2n<\sqrt{2}\,m$, there exists a $p$-Willmore Hopf tori with non-constant mean curvature in $\mathbb{S}^3$. These Hopf tori are based on the (non-trivial) closed $p$-elastic curves $\gamma_{n,m}$. Moreover, arguing as in \cite{Pi}, we obtain from the instability of closed $p$-elastic curves for $p\in(0,1)$ (Theorem \ref{instability}) that the corresponding $p$-Willmore Hopf tori are unstable as critical points of $\mathcal{W}_p$.

\begin{figure}[h!]
	\centering
	\includegraphics[height=5.55cm]{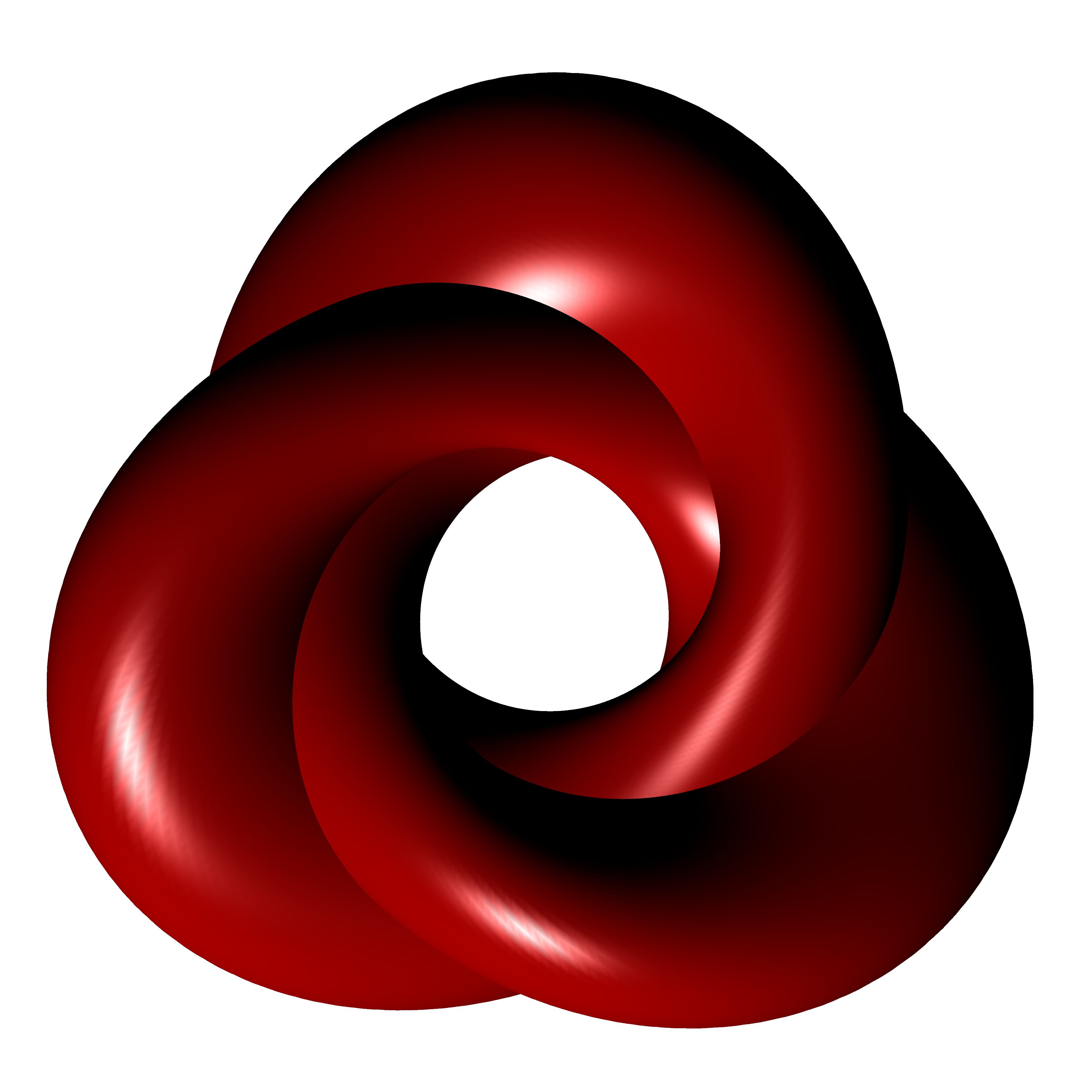}
	\includegraphics[height=5.55cm]{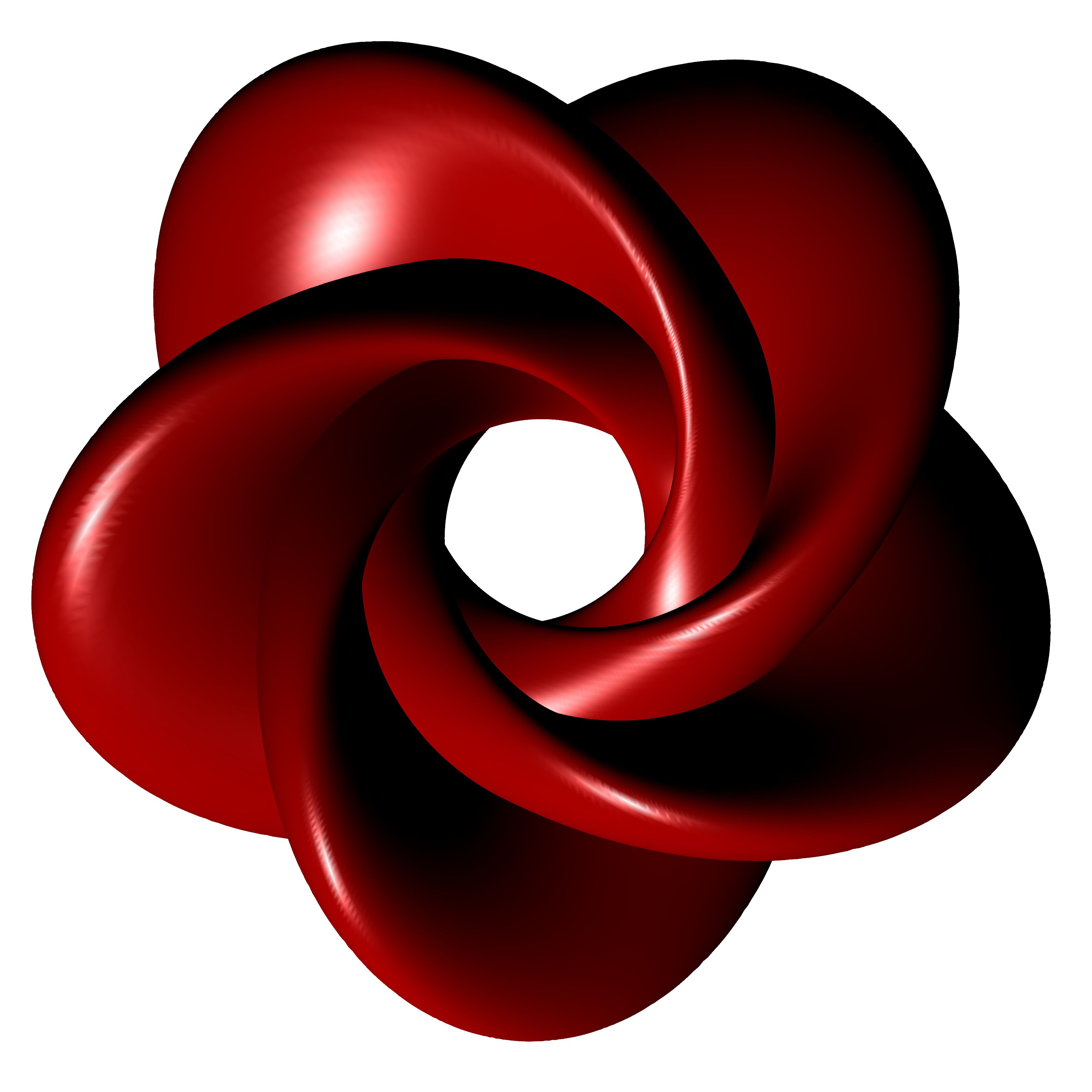}
	\includegraphics[height=5.55cm]{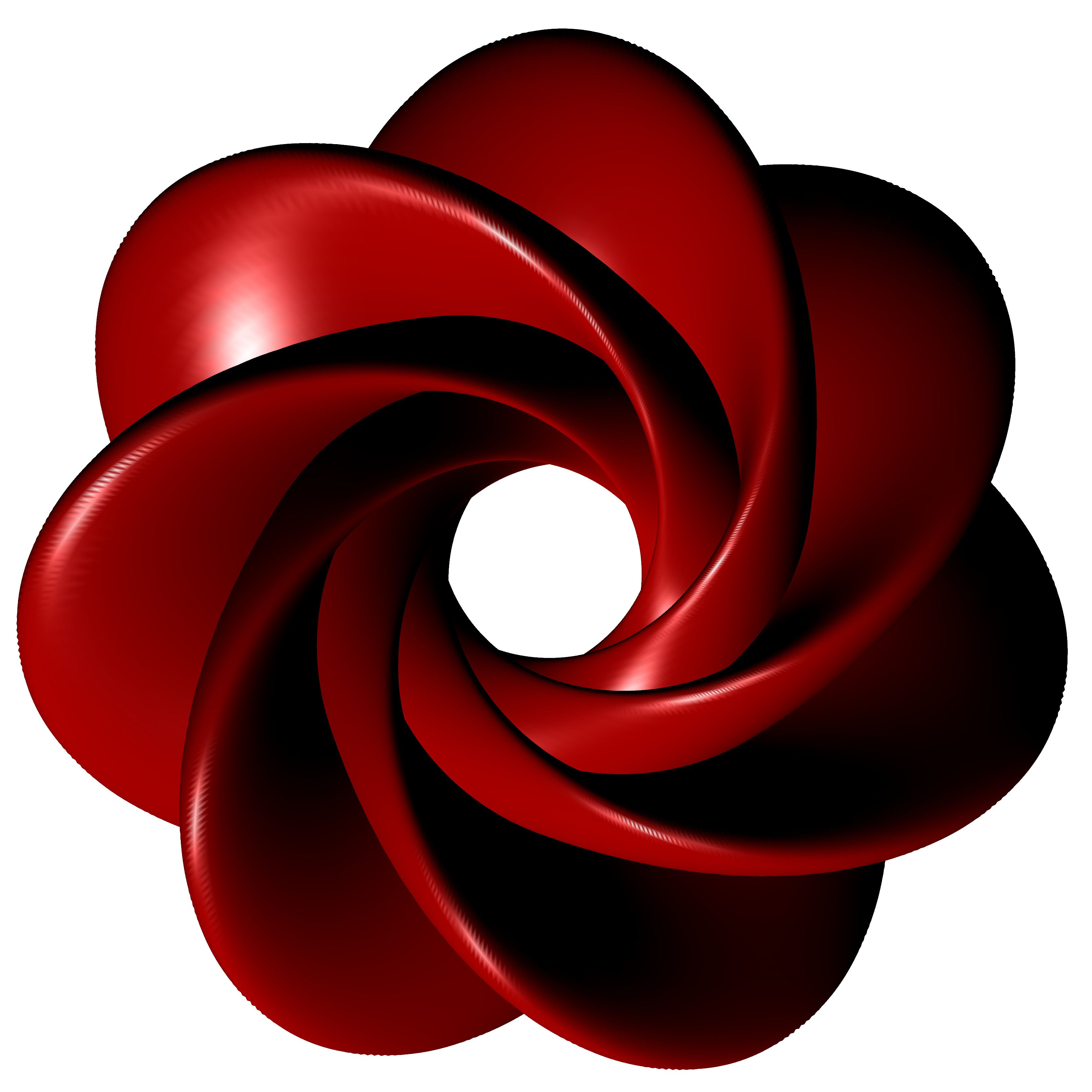}
	\caption{Stereogrpahic projection of the Hopf tori $X_\gamma$ whose base curves $\gamma$ are the $p$-elastic curves ($p=0.3$) in $\mathbb{S}^2$ represented in Figure \ref{F1}. The surfaces $X_\gamma$ are (unstable) $p$-Willmore tori for $p=0.3$ in $\mathbb{S}^3$.}
	\label{Hopf}
\end{figure}

In Figure \ref{Hopf}, using the parameterization of closed $p$-elastic curves given in \eqref{param} and of its horizontal lifts (see for instance \cite{P,PhD}), we illustrate the stereographic projection of three $p$-Willmore Hopf tori for $p=0.3$ whose base curves are represented in Figure \ref{F1}.

\section*{Acknowledgements}

We would like to thank Professor E. Aulisa for his valuable comments regarding the topic.

\begin{flushleft}
Anthony G{\footnotesize RUBER}\\
Department of Mathematics and Statistics, Texas Tech University, Lubbock, TX, 79409, USA\\
E-mail: anthony.gruber.d@gmail.com
\end{flushleft}

\begin{flushleft}
\'Alvaro P{\footnotesize \'AMPANO}\\
Department of Mathematics and Statistics, Texas Tech University, Lubbock, TX, 79409, USA\\
E-mail: alvaro.pampano@ttu.edu
\end{flushleft}

\begin{flushleft}
Magdalena T{\footnotesize ODA}\\
Department of Mathematics and Statistics, Texas Tech University, Lubbock, TX, 79409, USA\\
E-mail: magda.toda@ttu.edu
\end{flushleft}

\end{document}